\documentclass{elsart}
\usepackage{amssymb, amsmath, latexsym, amscd, graphics, epsfig}
\input xy
\input epsf

\xyoption{all}

\newcommand{\R}{\ensuremath{\mathbb{R}}}
\newcommand{\Z}{\ensuremath{\mathcal{Z}}}

\newcommand{\X}{\ensuremath{\mathbb{X}}}

\newcommand{\G}{\Gamma}

\newcommand{\cC}{\mathcal{C}}

\newcommand{\cS}{\mathcal{S}}
\newcommand{\cX}{\mathcal{X}}
\newcommand{\cY}{\mathcal{Y}}
\newcommand{\A}{\ensuremath{\mathcal{A}}}
\newcommand{\B}{\ensuremath{\mathcal{B}}}

\newcommand{\Ad}{\ensuremath{\operatorname{Ad}}}
\newcommand{\Vol}{\ensuremath{\operatorname{Vol}}}
\newcommand{\Aut}{\ensuremath{\operatorname{Aut}}}

\newcommand{\st}{\ensuremath{\operatorname{st}}}
\newcommand{\St}{\ensuremath{\operatorname{St}}}
\newcommand{\Ov}{\ensuremath{\operatorname{Over}}}
\newcommand{\Cov}{\ensuremath{\operatorname{Cov}}}

\newcommand{\bs}{\backslash}
\newcommand{\ba}{\backslash}
\newcommand{\quot}{\bs \! \bs}
\newcommand{\fgY}{\pi_1(G(Y))}
\newcommand{\fgYp}{\pi_1(G'(Y'))}
\newcommand{\fgYT}{\ensuremath{\pi_1(G(\cY),T)}}
\newcommand{\fgYpTp}{\ensuremath{\pi_1(G'(\cY'),T')}}
\newcommand{\fgYppTpp}{\ensuremath{\pi_1(G''(\cY''),T'')}}

\newcommand{\DY}{\widetilde{G(Y)}}
\newcommand{\DYp}{\widetilde{G'(Y')}}
\newcommand{\DYT}{D(\ensuremath{\cY,T)}}
\newcommand{\DYpTp}{D(\ensuremath{\cY',{T'})}}
\newcommand{\DYppTpp}{D(\ensuremath{\cY'',T'')}}

\newcommand{\Ovg}{\Ov(\G)}
\newcommand{\Covg}{\Cov(G(\Y))}

\def\Y{\mathcal{Y}}

\def\s{\sigma}
\def\G{\Gamma}
\def\X{\mathcal{X}}

\def\a{\alpha}

\def\s{\sigma}
\def\D{D(\mathcal{Y}, \phi)}

\begin{document}

\begin{frontmatter}

\title{Covering theory for complexes of groups} 
\author{Seonhee Lim}
\address{Mathematical Sciences Research Institute \\ 17 Gauss Way \\ Berkeley
CA 94720 \\ USA} \ead{seonheel@msri.org} 
\author{Anne Thomas\corauthref{cor}} \address{Mathematical
Sciences Research Institute \\ 17 Gauss Way \\ Berkeley CA 94720 \\ USA}
\ead{samohtenna@gmail.com}
\corauth[cor]{Corresponding author.  Telephone: +1 510 642 0610.  Fax: +1 510
642 0143.  From January 2008, contact details for Anne Thomas are:
Department of Mathematics, 310 Malott Hall, Cornell University, Ithaca NY
14853-4201, USA.  Telephone: +1 773 251 0972. Fax: +1 607 255 7149.}

\begin{abstract} We develop an explicit covering theory for
complexes of groups, parallel to that developed for graphs of groups
by Bass.  Given a covering of developable complexes of groups, we
construct the induced monomorphism of fundamental groups and
isometry of universal covers.  We characterize faithful complexes of
groups and prove a conjugacy theorem for groups acting freely on
polyhedral complexes.  We also define an equivalence relation on
coverings of complexes of groups, which allows us to construct a
bijection between such equivalence classes, and subgroups or
overgroups of a fixed lattice $\Gamma$ in the automorphism group of
a locally finite polyhedral complex $X$.
\end{abstract}


\end{frontmatter}

%****************************************************************************************************
%****************************************************************************************************
%****************************************************************************************************
\section{Introduction}
%****************************************************************************************************
%****************************************************************************************************
%****************************************************************************************************

Let $X$ be a locally finite polyhedral complex, such as a locally finite tree, 
Davis complex, or Bruhat--Tits building.  Then the group $G$ of automorphisms of
$X$ is naturally a locally compact group (see Section~\ref{ss:poly} below). 
Subgroups of $G$ with particular properties may be encoded by graphs of groups
(in the case of trees) or complexes of groups (for $\dim(X) \geq 2$) with
corresponding properties.  For example, $\G \leq G$ acting properly
discontinuously corresponds to a
graph or complex of finite groups, and $\G \leq G$ cocompact is encoded by a
finite complex of (possibly infinite) groups.  In this way, graphs or complexes
of groups may be used to study groups such as lattices in nonarchimedean Lie
groups or Kac-Moody groups, as well as lattices in automorphisms groups of products of trees, hyperbolic buildings, and so on.  A covering
theory for graphs of groups was developed by Bass in~\cite{b1:ctgg}.  This has
proved very useful for the study of tree lattices: see the reference
Bass--Lubotzky~\cite{BL}.  

The theory of complexes of groups is due to Gersten--Stallings~\cite{gs},
Corson~\cite{c} and Haefliger \cite{h1:cg},~\cite{BH}.  Haefliger, in Chapter
III.$\mathcal{C}$ of~\cite{BH}, translated into the framework of complexes of
groups the general theory of coverings of \'etale groupoids. (Chapter
III.$\mathcal{G}$ of~\cite{BH} discusses groupoids of local isometries.)  While
covering theory for \'etale groupoids is powerful, as it is strictly parallel to
the theory of coverings for topological spaces, the correspondence between
coverings of complexes of groups and coverings of \'etale groupoids is not
stated in~\cite{BH}.  Even the definition of the \'etale groupoid
canonically associated to a complex of groups is quite involved (pp. 595--596
of~\cite{BH}).  The \'etale groupoid perspective thus does not easily yield
results or constructions suitable for investigating concrete questions
concerning group actions on polyhedral complexes.  

One aim of this paper is to make
covering theory for complexes of groups more accessible, by following a more
explicit approach.  We also prove several results for group actions which
we hope will be broadly useful, including a
characterization of faithful complexes of groups (in
Section~\ref{ss:faithfulness}), and the Conjugacy Theorem
(Theorem~\ref{t:intro_conjugacy} below).  Finally, we establish in
Section~\ref{s:bijection} a bijection between suitably defined isomorphism
classes of coverings and subgroups or
overgroups of a fixed $\Gamma < \Aut(X)$.  This forms the technical background
for our work~\cite{LT} on counting overlattices, and hopefully will have other
applications.

Let us briefly recall Haefliger's theory of complexes of groups (see
Section~\ref{ss:cxs_of_gps} below for details, and in particular for the
definition of covering).  The action of a group $G$ on
a simply connected polyhedral complex $X$ induces a complex of groups $G(Y)$
over the quotient $Y = G \bs X$.  The fundamental group $\pi_1(G(Y))$  then
acts on the simply connected universal cover $\DY$, with $\fgY$ isomorphic to
$G$, and $\DY$ isometric to $X$.  An arbitrary complex of groups $G(Y)$ is
developable if it is induced by a group action in this way. A key difference
between Bass--Serre theory and the theory of complexes of groups is that complexes of groups
need not be developable.  However, if a complex of groups has nonpositive
curvature (see Section~\ref{sss:nonpos}), it is developable.

Our first main result describes the functoriality of coverings.

\begin{thm}\label{t:functor_cov}
 Let $\lambda:G(Y)\to G'(Y')$ be a covering of developable complexes of
 groups.  Then $\lambda$ induces a monomorphism of fundamental groups
\[\Lambda: \fgY \rightarrow \fgYp\] and a $\Lambda$-equivariant isometry of
universal covers
\[L:\DY \rightarrow \DYp. \]
\end{thm}

\noindent Theorem~\ref{t:functor_cov} also follows from covering theory for
\'etale groupoids (Haefliger, personal communication).  Our contribution is to construct the maps $\Lambda$ and $L$
explicitly; we then make repeated use of these constructions in later sections of this
work.  Theorem~\ref{t:functor_cov} is proved in Section~\ref{ss:functor_cov},
using material from Section~\ref{ss:functor_morphism}.

In Section~\ref{ss:faithfulness}, we characterize the group
\[N = \ker\left(\fgY \rightarrow \Aut(\DY)\right)\] where $G(Y)$ is developable.
If $N$ is trivial, then the
complex of groups $G(Y)$ is said to be \emph{faithful}, and we may
identify the fundamental group $\pi_1(G(Y))$ with a subgroup of
$\Aut(\DY)$.

In Section~\ref{ss:main_lemma} we develop technical results, similar to those in
Section~4 of~\cite{b1:ctgg} in the case of trees. As
described in Proposition~2.1 of~\cite{T}, Haefliger's morphisms of complexes of
groups,
when restricted to complexes of groups over $1$--dimensional spaces, are not the
same as Bass' morphisms of graphs of groups.  Also, the universal covers of graphs of
groups and of complexes of groups are defined with respect to different choices.
Hence, our proofs differ in many details from those of~\cite{b1:ctgg}.

An additional consideration for complexes of groups, which has no
analogue in Bass--Serre theory, is the relationship between
coverings and developability.  In
Section~\ref{ss:cov_dev}, we show:

\begin{prop}\label{p:covdev}  Let $\lambda:G(Y)\to G'(Y')$ be a covering of complexes of
groups.
\begin{enumerate} \item
If $G'(Y')$ is developable, then $G(Y)$ is developable.
\item If $G(Y)$ has nonpositive curvature (hence is developable),
then $G'(Y')$ has nonpositive curvature, hence $G'(Y')$ is
developable.
\end{enumerate}
\end{prop}

One of the main applications of the results of Section~\ref{ss:main_lemma} is
the Conjugacy Theorem below, proved as
Theorem~\ref{t:conjugacy} in Section~\ref{s:conjugacy}.  Let $H$ be a subgroup
(acting without inversions) of $G=\Aut(X)$, for $X$ a locally finite polyhedral
complex, and define \[G_H = \{ g \in G \mid g\sigma \in H\sigma \mbox{ for all
cells $\sigma$ of $X$} \}.\]

\begin{thm}[Conjugacy Theorem]\label{t:intro_conjugacy} If $\G \leq G_H$ acts freely on $X$ then there is an element
$g \in G_H$ such that $g \G g^{-1} \leq H$.
\end{thm}

\noindent The corresponding result for trees (Theorem~5.2 of~\cite{b1:ctgg}) was a basic
tool in~\cite{BK}.  In~\cite{b1:ctgg}, as well as a proof using covering theory
for graphs of groups, a simple direct proof due to the referee was given.  This
relied on the fact that a group acting freely on a tree is free.  In higher
dimensions, it seems that covering theory must be used.

In Section~\ref{s:bijection} we define isomorphism of coverings (see
Definition~\ref{d:isom_covering}) so that the following bijection holds:

\begin{thm}\label{t:bijection} Let $X$ be a simply connected polyhedral
complex, and let $\G$ be a subgroup of $\Aut(X)$ (acting without
inversions) which induces a complex of groups $G(Y)$. Then there is
a bijection between the set of subgroups of $\Aut(X)$ (acting
without inversions) which contain $\G$, and the set of isomorphism
classes of coverings of faithful, developable complexes of groups by
$G(Y)$.
\end{thm}

\noindent The main ingredients in the proof of Theorem~\ref{t:bijection} are
Theorem~\ref{t:functor_cov} above, and the results of
Section~\ref{ss:main_lemma}.  As a corollary to Theorem~\ref{t:bijection}, we
show that there is a bijection between $n$--sheeted coverings, and overlattices
of index $n$ (that is, lattices containing a fixed lattice $\G$ with index
$n$).  Similar results hold for subgroups and sublattices.  In~\cite{L1}, Lim
defined isomorphism of coverings of graphs of groups and proved the bijection of
Theorem~\ref{t:bijection} for trees.

%****************************************************************************************************
%****************************************************************************************************
%****************************************************************************************************
\section{Background}\label{s:background}
%****************************************************************************************************
%****************************************************************************************************
%****************************************************************************************************

We begin by recalling the basic theory of lattices, in
Section~\ref{ss:lattices}. Since the quotient of a simplicial
complex by a simplicial group action is not in general a simplicial
complex, it is natural to define complexes of groups over polyhedral
complexes instead.  In Section~\ref{ss:poly}, we give definitions of polyhedral
complexes and the topology of their automorphism groups. Small
categories without loops, or scwols, are algebraic objects that
substitute for polyhedral complexes. These are described in
Section~\ref{ss:scwols} (following section III.$\cC$ 1-2 of~\cite{BH}).
 We conclude this background material by, in
Section~\ref{ss:cxs_of_gps}, summarizing Haefliger's theory of
complexes of groups, as presented in Chapter~III.$\cC$ of~\cite{BH}.

%****************************************************************************************************
%****************************************************************************************************
\subsection{Lattices}\label{ss:lattices}
%****************************************************************************************************
%****************************************************************************************************

Let $G$ be a locally compact topological group with left-invariant
Haar measure $\mu$. A discrete subgroup $\Gamma$ of $G$ is a
\emph{lattice} if its covolume $\mu(\Gamma \bs G)$ is finite. A
lattice is called \textit{cocompact} or \emph{uniform} if $\Gamma
\backslash G$ is compact.

Let $\mathcal{S}$ be a left $G$-set such that, for each $s \in
\mathcal{S}$, the stabilizer $G_s$ is compact and open. For any
discrete subgroup $\Gamma$ of $G$, the stabilizers $\G_s$ are finite
groups, and we define the \emph{$\mathcal{S}$-covolume} of $\G$ as
$$\Vol(\Gamma \quot \mathcal{S})= \sum_{s \in \Gamma \bs \mathcal{S}}
\frac{1}{|\Gamma_s|} \leq \infty.$$ It is shown in \cite{BL},
Chapter~1, that if $G\bs \cS$ is finite and $G$ admits a lattice,
then there is a normalization of the Haar measure $\mu$, depending
only on $\mathcal{S}$, such that for every discrete subgroup
$\Gamma$ of $G$, \[ \mu(\Gamma \backslash G)=\Vol(\Gamma \quot
\mathcal{S}).\]  It is clear that for two lattices $\Gamma \subset
\Gamma'$ of $G$, the index $[\Gamma' : \Gamma]$ is equal to the
ratio of the covolumes $\mu(\Gamma \bs G):\mu(\Gamma' \bs G)$.

%****************************************************************************************************
%****************************************************************************************************
\subsection{Polyhedral complexes}\label{ss:poly}
%****************************************************************************************************
%****************************************************************************************************

Let $M_\kappa^n$ be the complete, simply connected, Riemannian
$n$-manifold of constant sectional curvature $\kappa \in \R$.

\begin{defn}[polyhedral complex] An $M_\kappa$-polyhedral complex $K$ is a finite-dimensional CW-complex such that: \begin{enumerate}
\item each open cell of dimension $n$ is isometric to the interior of a compact convex polyhedron in $M^n_\kappa$;
and \item for each cell $\s$ of $K$, the restriction of the
attaching map to each open codimension one face of $\s$ is an
isometry onto an open cell of $K$.\end{enumerate}\end{defn}

If an $M_\kappa$-polyhedral complex is locally finite, then it is
a geodesic metric space by the Hopf--Rinow Theorem (see, for
example,~\cite{BH}).  More generally, we have:

\begin{thm}[Bridson,~\cite{BH}] An
$M_\kappa$-polyhedral complex with finitely many isometry classes of cells is a
complete geodesic metric space.
\end{thm}

Let $K$ be a locally finite, connected polyhedral complex, and let
$\Aut(K)$ be the group of cellular isometries, or automorphisms, of
$K$.  Then $\Aut(K)$ is naturally a locally compact group, with a
neighborhood basis of the identity consisting of automorphisms
fixing larger and larger balls.  With respect to this topology, a
subgroup $\G$ of $\Aut(K)$ is discrete if and only if for each cell
$\s$ of $K$, the stabilizer $\G_\s$ is finite.   A subgroup $\G$ of
$\Aut(K)$ is said to act \emph{without
 inversions} if whenever $g \in  \G$ preserves a cell of $K$, $g$ fixes that cell pointwise.

%****************************************************************************************************
%****************************************************************************************************
\subsection{Small categories without loops}\label{ss:scwols}
%****************************************************************************************************
%****************************************************************************************************

In Chapter~III.$\cC$ of~\cite{BH}, complexes of
groups are presented using the language of scwols, or small categories without
loops.  As we explain in this section, to any polyhedral complex $K$
one may associate a scwol $\X$, which has a geometric realization
$|\X|$ isometric to the barycentric subdivision of $K$.  Morphisms
of scwols correspond to polyhedral maps, and group actions on scwols
correspond to actions without inversions on polyhedral complexes.

\begin{defn}[scwol] A small category without loops (scwol) $\X$ is a disjoint union of a set
$V(\X)$, the vertex set, and a set $E(\X)$, the edge set, endowed
with maps \[i: E(\X) \to V(\X) \quad\mbox{and}\quad t: E(\X) \to
V(\X)\] and, if $E^{(2)}(\X)$ denotes the set of pairs $(a,b)$ of
edges where
$i(a)=t(b)$, with a map \begin{align*}E^{(2)}(\X) & \to E(\X) \\
(a,b) & \mapsto ab\end{align*} such that:
\begin{enumerate}
\item if $(a,b) \in E^{(2)}(\X)$, then $i(ab)=i(b)$ and
$t(ab)=t(a)$; \item if $a,b,c \in E(\X)$ such that $i(a)=t(b)$ and
$i(b)=t(c)$, then $(ab)c=a(bc)$; and \item if $a \in E(\X)$, then
$i(a) \neq t(a)$.\end{enumerate}
\end{defn}

For $a \in E(\X)$, the vertices $i(a)$ and $t(a)$ are called the
\emph{initial vertex} and \emph{terminal vertex} of $a$
respectively.  If $(a,b) \in E^{(2)}(\X)$ we say $a$ and $b$ are
\emph{composable}, and that $ab$ is the \emph{composition} of $a$
and $b$. We will sometimes write $\alpha \in \X$ for $\alpha \in
V(\X) \cup E(\X)$. If $\alpha \in V(\X)$ then $i(\alpha)=t(\alpha) =
\alpha$.

The motivating example of a scwol is the scwol $\X$ associated to a
polyhedral complex $K$. The set of vertices $V(\X)$ corresponds to
the set of cells of $K$ (or the set of barycenters of the cells of
$K$). The set of edges $E(\X)$ is the set of $1$-simplices of the
barycentric subdivision of $K$, that is, each element of $E(\X)$
corresponds to a pair of cells $T \subsetneq S$, with initial vertex
$S$ and terminal vertex $T$.  The composition of the edge $a$
corresponding to $T \subsetneq S$ and the edge $b$ corresponding to
$S \subsetneq U$ is the edge $ab$ corresponding to $T \subsetneq U$.

\begin{figure}[ht]
\begin{center}
\includegraphics{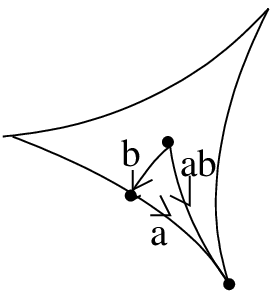}
\end{center}
\end{figure}

Conversely, given a scwol $\X$, we may construct a polyhedral
complex, called the \emph{geometric realization}. For an integer $k
\geq 0$, let $E^{(k)}(\X)$ be the set of sequences $(a_1, a_2,
\ldots, a_k)$ of composable edges, that is, $(a_{j}, a_{j+1}) \in
E^{(2)}(\X)$ if $k > 1$, $E^{(1)}(\X) = E(\X)$, and $E^{(0)}(\X) =
V(\X)$. The geometric realization $|\X|$ of $\X$ is defined as a
polyhedral complex whose cells of dimension $k$ are standard
$k$-simplices indexed by the elements of $E^{(k)}(\X)$. For the
details of this construction, see~\cite{BH}, pp. 522--523. If $\X$
is the scwol associated to an $M_\kappa$-polyhedral complex $K$,
then $|\X|$ may be realized as an $M_\kappa$-polyhedral complex
isometric to the barycentric subdivision of $K$.

For a scwol $\X$, let $E^\pm(\X)$ be the set of \emph{oriented
edges}, that is, the set of symbols $a^+$ and $a^-$, where $a \in
E(\X)$. For $e=a^+$, we define $i(e)=t(a)$, $t(e)=i(a)$ and
$e^{-1}=a^-$. For $e=a^-$, we define $i(e)=i(a)$, $t(e)=t(a)$ and
$e^{-1}=a^+$.

An \emph{edge path} in $\X$ joining the vertex $\s$ to the vertex
$\tau$ is a sequence $(e_1, e_2, \ldots, e_n)$ of elements of $E^\pm
(\X)$ such that $i(e_1)=\s$, $i(e_{j+1})=t(e_j)$ for $1 \leq j \leq
n-1$ and $t(e_n)=\tau$.

A scwol $\X$ is \emph{connected} if for any two vertices $\s, \tau
\in V(\X)$, there is an edge path joining $\s$ to $\tau$.
Equivalently, $\X$ is connected if and only if the geometric
realization $|\X|$ is connected.  A scwol is \emph{simply connected}
if and only if its geometric realization is simply connected as a
topological space.

\begin{defn}[morphism of scwols]\label{d:morphism_scwols} Let $\X$ and $\X'$ be two scwols. A morphism $l: \X \to \X'$ is
a  map that sends $V(\X)$ to $V(\X')$ and $E(\X)$ to $E(\X')$, such
that
\begin{enumerate} \item\label{i:adj} for each $a \in E(\X)$, we have
$i(l(a))=l(i(a))$ and $t(l(a))=l(t(a))$; and
\item\label{i:comp} for each $(a,b) \in E^{(2)}(\X)$, we have $l(ab)=l(a)l(b)$.
\newline \noindent A nondegenerate morphism of scwols is a morphism of scwols such that
in addition to~\eqref{i:adj} and~\eqref{i:comp},
\item\label{i:bij} for each vertex $\sigma \in V(\X)$, the restriction of $l$ to
the set of edges with initial vertex $\sigma$ is a bijection onto
the set of edges of $\X'$ with initial vertex
$l(\sigma)$.\end{enumerate}
\end{defn}

An \emph{automorphism} of a scwol $\X$ is a morphism $l:\X \to \X$
which has an inverse.  Note that Condition~\eqref{i:bij} in
Definition~\ref{d:morphism_scwols} is automatic for automorphisms.

\begin{defn}[covering of scwols]\label{d:cov_scwols} Let $\X$ be a (nonempty) scwol and let $\X'$ be a connected
scwol.  A nondegenerate morphism of scwols $l:\X \to \X'$ is called
a covering if, for every vertex $\sigma$ of $\X$, the restriction of
$l$ to the set of edges with terminal vertex $\sigma$ is a bijection
onto the set of edges of $\X'$ with terminal vertex $l(\sigma)$.
\end{defn}

Let $\X$ and $\X'$ be scwols associated to polyhedral complexes $K$
and $K'$ respectively.  A polyhedral map $K \to K'$ induces a
morphism of scwols $\X \to \X'$, and conversely, a morphism $l: \X
\to \X'$ induces a continuous polyhedral map $|l|:|\X|\to|\X'|$ (see
\cite{BH}, p. 526).  The morphism $l$ is nondegenerate if and only
if the restriction of $|l|$ to the interior of each cell of $K$
induces a homeomorphism onto the interior of a cell of $K'$, and $l$
is a covering if and only if $|l|$ is a (topological) covering. A
morphism $l:\X \to \X$ is an automorphism of $\X$ if and only if
$|l|:K \to K$ is an automorphism of $K$.

\begin{defn}[group actions on scwols]\label{d:actiononscwol} An action of a group $G$ on a scwol $\X$ is a homomorphism from
$G$ to the group of automorphisms of $\X$ such
that:\begin{enumerate}\item for all $a \in E(\X)$ and $g \in G$, we
have $g \cdot i(a) \neq t(a)$; and \item for all $g \in G$ and $a
\in E(\X)$, if $g \cdot i(a)=i(a)$ then $g \cdot a=a$ (no
``inversions").
\end{enumerate}
\end{defn}

The action of a group $G$ on a scwol $\X$ induces a \emph{quotient
scwol} $\Y = G \bs \X$, defined as follows.  The vertex set is
$V(\Y) = G \bs V(\X)$ and the edge set $E(\Y) = G\bs E(\X)$.  For
every $a \in E(\X)$ we have $i(Ga) = Gi(a)$ and $t(Ga) = Gt(a)$, and
if $(a,b) \in E^{(2)}(\X)$ then the composition of $Ga$ and $Gb$ is
$Gab$.  The natural projection $p:\X \to \Y$ is a nondegenerate
morphism of scwols.

Let $\X$ be the scwol associated to a polyhedral complex $K$, and
let $\G$ be a subgroup of $G=\Aut(K)$. Then $\G$ acts on $\X$, in
the sense of Definition~\ref{d:actiononscwol}, if and only if $\G$
acts without inversions on $K$.

In the case $K$ is locally finite, we define the covolume of a
discrete subgroup $\G\leq G$ acting on $\X$ as follows. For the
$\G$-set $\mathcal{S}$ in Section~\ref{ss:lattices}, we choose the
set of vertices $V(\X)$ (which corresponds to the set of cells of
$K$).
 By the same
arguments as for tree lattices (\cite{BL}, Chapter~1), it can be
shown that if $G \bs K$ is finite, then $\G$ is a lattice if and
only if its $V(\X)$-covolume converges, and $\G$ is a cocompact
lattice if and only if $\G \bs V(\X)$ is a finite set.  We now
normalize the Haar measure $\mu$ on $G$ so that
$$\mu(\Gamma \backslash G) =\Vol(\G \quot V(\X)) = \sum_{\s \in \Gamma \backslash
V(\X)} \frac{1}{|\Gamma_\s|}.$$

%****************************************************************************************************
%****************************************************************************************************
\subsection{Complexes of groups}\label{ss:cxs_of_gps}
%****************************************************************************************************
%****************************************************************************************************

In this section, we recall Haefliger's theory of complexes of
groups.  We mainly follow the notation and definitions of Chapter
III.$\mathcal{C}$ of~\cite{BH}, although at times, such as in Defintion~\ref{d:induced_morphism} and
Proposition~\ref{p:isoms} below, we indicate choices and define maps
more explicitly.  Section~\ref{sss:objects_morphisms}  defines
complexes of groups and their morphisms.  Section~\ref{sss:fundgp}
then discusses groups associated to complexes of groups, in
particular the fundamental group, and Section~\ref{sss:univ_cover}
discusses scwols associated to complexes of groups, in particular
the universal cover. In Section~\ref{sss:nonpos} we describe the
role of local developments and nonpositive curvature. All references
to~\cite{BH} in this section are to Chapter III.$\mathcal{C}$, which
the reader should consult for further details.

%*****************************************************************************
\subsubsection{Objects and morphisms of the category of complexes of
groups}\label{sss:objects_morphisms}
%********************************************************************************

\begin{defn}[complex of groups]\label{d:cx_of_gps} Let $\Y$ be a scwol. A complex of groups $G(\Y)=(G_\s, \psi_a, g_{a,b})$
over $\Y$ is given by the following data:
\begin{enumerate}
\item for each $\s \in V(\Y)$, a group $G_\s$, called the local
group at $\s$;
\item for each $a \in E(\Y)$, an injective group homomorphism $\psi_a : G_{i(a)} \to
G_{t(a)}$; and
\item for each pair of composable edges $(a,b) \in E^{(2)}(\Y)$,
a twisting element $g_{a,b} \in G_{t(a)}$;
\end{enumerate}
with the following properties:
\renewcommand{\labelenumi}{{\normalfont (\roman{enumi})}}
\begin{enumerate} \item $Ad(g_{a,b})
\psi_{ab} =\psi_a \psi_b$, where $Ad(g_{a,b})$ denotes conjugation
by $g_{a,b}$; and
\item $\psi_a (g_{b,c})g_{a, bc} = g_{a,b} g_{ab,c}$, for each triple
$(a,b,c) \in E^{(3)}(\Y)$.
\end{enumerate}
\end{defn}

For example, any group $G$ is a complex of groups over a singleton
$\Y = \{ \ast \} = V(\Y)$, with $G_\ast = G$; since $E(\Y) = \phi$,
no other data is necessary.

\begin{defn}[morphism of complexes of groups]\label{d:morphism_cx_of_gps} Let $G(\Y)$ be as in Definition~\ref{d:cx_of_gps}
 and let $G'(\Y')=(G'_{\s'}, \psi_{a'}, g_{a',b'})$ be another
 complex of groups over a scwol $\Y'$. Let $l: \Y \to \Y'$ be a morphism of scwols. A morphism $\phi=(\phi_\s,
 \phi(a)): G(\Y) \to G'(\Y')$ of complexes of groups over $l$ consists of
\begin{enumerate}
\item a group homomorphism $\phi_\s : G_\s \to G'_{l(\s)}$, called
the local map at $\s$, for each $\s \in V(\Y)$; and \item an element
$\phi(a) \in G'_{t(l(a))}$ for each $a \in E(\Y)$; \end{enumerate}
such that:\renewcommand{\labelenumi}{{\normalfont (\roman{enumi})}}
\begin{enumerate}
\item $Ad(\phi(a))\psi_{l(a)}\phi_{i(a)}=\phi_{t(a)} \psi_a$; and
\item $\phi_{t(a)}(g_{a,b})\phi(ab)=\phi(a)
\psi_{l(a)}(\phi(b))g_{l(a), l(b),}$ for every $(a,b) \in
E^{(2)}(\Y)$.\end{enumerate}
\end{defn}

A morphism $\phi$ is an \emph{isomorphism} if $l$ is an isomorphism
of scwols and $\phi_\s$ is a group isomorphism for every $\s \in
V(\Y)$.  A morphism $\phi$ is \emph{injective on the local groups}
if each of the maps $\phi_\s$ is injective.

The \emph{composition} $\phi' \circ \phi$ of a morphism
$\phi=(\phi_\s, \phi(a)):G(\Y) \to G'(\Y')$ over $l$ and a morphism
$\phi'=(\phi'_\s, \phi'(a)):G'(\Y') \to G''(\Y'')$ over $l'$ is the
morphism over $l'\circ l$ defined by the homomorphisms $(\phi' \circ
\phi)_\s = \phi'_{l(\s)} \circ \phi_\s$ and the elements
$(\phi'\circ \phi)(a) = \phi'_{l(t(a))}(\phi(a))\phi'(l(a))$.

A special case of a morphism of complexes of groups is when $\Y'$ in
Definition~\ref{d:morphism_cx_of_gps} is a singleton, with $G'_\ast
= G'$. In this case, $\phi$ may be regarded as a morphism from the
complex of groups $G(\Y)$ to the group $G'$.

\begin{defn}[homotopy]\label{d:homotopy} Let $\phi$ and $\phi'$ be
two morphisms from $G(\Y)$ to a group $G'$, given respectively by
$(\phi_\s, \phi(a))$ and $(\phi'_\s, \phi'(a))$.  A homotopy from
$\phi$ to $\phi'$ is given by a family of elements $k_\s \in G'$,
indexed by $\s \in V(\Y)$, such that
\begin{enumerate}
\item $\phi'_\s = \Ad(k_\s)\phi_\s$ for all $\s \in V(\Y)$; and
\item $\phi'(a) = k_{t(a)} \phi(a) k_{i(a)}^{-1}$ for
all $a \in E(\Y)$.
\end{enumerate}
\end{defn}

Let $G$ be a group acting on a scwol $\X$ with quotient $\Y = G \bs
\X$, and let $p: \X \to \Y$ be the natural projection. The
\emph{complex of groups $G(\Y)=(G_\s, \psi_a, g_{a,b})$ associated
to the action of $G$ on $\X$} is defined as follows.

For each vertex $\s \in V(\Y)$, choose a vertex $\overline{\s} \in
V(\X)$ such that $p(\overline\s) =\s$.  For each edge $a \in E(\Y)$
with $i(a)=\s$, there exists a unique edge $\overline{a} \in E(\X)$
such that $p(\overline{a})=a$ and $i(\overline{a})=\overline\s$.
Choose $h_a \in G$ such that $h_a \cdot
t(\overline{a})=\overline{t(a)}$. For each $\sigma \in V(\Y)$, let
$G_\s$ be the stabilizer in $G$ of $\overline{\s} \in V(\X)$.  For
each $a \in E(\Y)$, let $\psi_a: G_{i(a)} \to G_{t(a)}$ be
conjugation by $h_a$, that is, \[ \psi_a: g \mapsto h_a g h_a
^{-1}.\]  For every pair of composable edges $(a,b) \in E^{(2)}(\Y)$,
define $g_{a,b}=h_a h_b h_{ab} ^{-1}$. Then
$G(\Y)=(G_\s,\psi_a,g_{a,b})$ is a complex of groups.

When precision is needed, we denote the set of choices of
$\overline\s$ and $h_a$ in this construction by $C_\bullet$, and the
complex of groups $G(\Y)$ constructed with respect to these choices
by $G(\Y)_{C_\bullet}$.  If $C_\bullet'$ is another choice of
$\overline\s'$, $h_a'$, then an isomorphism $\phi =
(\phi_\s,\phi(a))$ from $G(\Y)_{C_\bullet}$ to $G(\Y)_{C_\bullet'}$
is obtained by choosing elements $k_\s \in G$, such that for each
$\s \in V(\Y)$, $k_\s \cdot \overline\s = \overline\s'$.  Then put
$\phi_\s = \Ad(k_\s)|_{G_\s}$ and $\phi(a) = k_{t(a)}h_a
k_{i(a)}^{-1}h_a'^{-1}$.

When $G(\Y)$ is a complex of groups associated to an action of a
group $G$, there is a canonical morphism of complexes of groups
$\phi_1 : G(\Y) \to G$, given by $\phi_1=(\phi_\sigma, \phi(a))$,
with $\phi_\sigma=G_\sigma \to G$ the inclusion, and $\phi(a)=h_a$.

\begin{defn}[developable] A complex of groups $G(\Y)$ is developable if it is
isomorphic to a complex of groups associated to the action of a
group $G$ on a scwol $\X$ in the above sense, with $\Y = G \bs
\X$.\end{defn}

\begin{prop}[Corollary~2.15,~\cite{BH}]\label{p:dev_inj}
A complex of groups $G(\Y)$ is developable if and only if there
exists a morphism $\phi$ from $G(\Y)$ to some group $G$ which is
injective on the local groups.\end{prop}

We now define coverings.

\begin{defn}[covering of complexes of groups]\label{d:covering}
Let $\phi:G(\Y) \to G'(\Y')$ be a morphism of complexes of groups
over a nondegenerate morphism of scwols $l: \Y \to \Y'$, where $\Y'$
is connected. The morphism $\phi$ is a covering (of $G'(\Y')$ by
$G(\Y)$) if for each vertex $\s \in V(\Y)$,
\begin{enumerate}
\item  the group homomorphism $\phi_{\s}: G_{\s} \to G'_{l(\s)}$ is
injective, and \item\label{i:cov_bijection} for every $a' \in
E(\Y')$ and $\s \in V(\Y)$ with $t(a')=\s'=l(\s)$, the map
\[
\coprod_{\substack{a \in l^{-1}(a') \\ t(a)=\s}} G_\s
/\psi_{a}(G_{i(a)}) \to G'_{\s'}/\psi_{a'} (G'_{i(a')})\] induced by
\[
 g \mapsto \phi_{\s}(g) \phi(a)
 \] is bijective.
\end{enumerate}\end{defn}

\noindent{}From Condition~\eqref{i:cov_bijection} of this
definition, it follows that
$$ \sum_{\substack{ a \in  l^{-1}(a') \\ t(a) =\s}} \frac{|G_{\s }|}{|G_{i(a)}|}= \frac{|G'_{\s'}|}{|
G'_{i(a')}|}.$$ Since $\Y'$ is connected, the value of
$$n:= \sum_{ \s \in \; l^{-1}(\s')} \frac{|G'_{\s'}|}{|G_\s|}= 
\sum_{ a \in l^{-1}(a')}
\frac{|G'_{i(a')}|}{|G_{i(a)}|}$$ is independent of the vertex $\s'$
and the edge $a'$.  A covering of complexes of groups with the above
$n$ is said to be \emph{$n$--sheeted}.

We will often use Definition~\ref{d:induced_morphism} below, which
defines a morphism of complexes of groups induced by an equivariant
morphism of scwols, keeping track of the choices we make.

\begin{defn}[induced morphism]\label{d:induced_morphism} Let $\X$ and $\X'$ be simply connected scwols,
endowed with actions of groups $G$ and $G'$, and let $\Y = G \bs \X$
and $\Y'=G'\bs \X'$ be the quotient scwols.  Let $L: \mathcal{X} \to
\mathcal{X}'$ be a morphism of scwols which is equivariant with
respect to a group morphism $\Lambda: G \to G'$. Let $l: \mathcal{Y}
\to \Y'$ be the induced morphism of the quotients.

For any choices $C_{\bullet}$ and $C'_{\bullet}$ of data for the
actions of $G$ and $G'$ on $\X$ and $\X'$, and for any choice
$N_{\bullet}$ of elements $k_{\sigma} \in G'$ indexed by $\sigma \in
V(\mathcal{Y})$ such that $k_{\sigma}\cdot
L(\overline{\sigma})=\overline{l(\sigma)}$, there is an associated
morphism of complexes of groups
\[\lambda=\lambda_{C_{\bullet}, C'_{\bullet}, N_{\bullet}}:
G(\mathcal{Y})_{C_{\bullet}} \to G'(\mathcal{Y}')_{C'_{\bullet}}\]
over $l$, given by
\begin{align*}
\lambda_{\sigma}: G_{\sigma} &\to G'_{l(\sigma)}\\
g &\mapsto k_{\sigma}\Lambda(g) k_{\sigma}^{-1}
\end{align*}
and
\begin{equation*}
\lambda(a)=k_{t(a)} \Lambda(h_a) k^{-1}_{i(a)}h'^{-1}_{l(a)}
\end{equation*}
(see Section~2.9(4),~\cite{BH}). \end{defn}

%*********************************************************************************
\subsubsection{The fundamental group of a complex of
groups}\label{sss:fundgp}
%*********************************************************************************

There are two definitions of the fundamental group of a complex of
groups, which result in canonically isomorphic groups.  Both
definitions involve the universal group.

\begin{defn}[universal group]\label{d:univ_gp} The universal group $FG(\Y)$ of a complex of groups $G(\Y)$ over a scwol $\Y$ is the group presented by the generating set
\[\coprod_{ \s \in V(\Y)} G_\s \coprod E^{\pm}(\Y)\] with the
following relations:
\begin{enumerate}
\item the relations in the groups $G_\s$;
\item $(a^+)^{-1}=a^-$ and $(a^-)^{-1}=a^+$;
\item $a^+ b^+ = g_{a,b}(ab)^+$, for every $(a,b) \in E^{(2)}(\Y)$;
and \item $\psi_a (g) = a^+ g a^-$, for every $g \in G_{i(a)}$.
\end{enumerate}
\end{defn}

There is a natural morphism $\iota=(\iota_\s, \iota(a)): G(\Y) \to
FG(\Y)$, where $\iota_\s: G_\s \to FG(\Y)$ takes $G_\s$ to its image
in $FG(\Y)$, and $\iota(a)=a^+$.

\begin{prop}[Proposition~3.9,~\cite{BH}]\label{p:dev_inj_FG}
A complex of groups $G(\Y)$ over a connected scwol $\Y$ is
developable if and only if $\iota:G(\Y) \to FG(\Y)$ is injective on
the local groups.\end{prop}

The first definition of the fundamental group of a complex of groups
$G(\Y)$ involves the choice of a basepoint $\s_0 \in V(\Y)$. A
\emph{$G(\Y)$-path} starting from $\s_0$ is then a sequence $(g_0,
e_1, g_1, e_2, \ldots, e_n, g_n)$ where $(e_1, e_2, \ldots, e_n)$ is
an edge path in $\Y$ starting from $\s_0$, we have $g_0 \in
G_{\s_0}$, and $g_j \in G_{t(e_j)}$ for $1 \leq j \leq n$. A
$G(\Y)$-path joining $\s_0$ to $\s_0$ is called a
\emph{$G(\Y)$--loop} at $\s_0$.

To each path $c=(g_0, e_1, g_1, e_2, \ldots, e_n, g_n)$, we
associate the element $\pi(c)$ of $FG(\Y)$ represented by the word
$g_0 e_1 g_1 \cdots e_n g_n$.  Suppose now that $c$  and $c'=(g'_0,
e'_1, g'_1, \ldots, e'_n, g'_n)$ are two $G(\Y)$--loops at $\s_0$. We
say $c$ and $c'$ are \emph{homotopic} if $\pi(c)=\pi(c')$, and
denote the homotopy class of $c$ by $[c]$. The \emph{concatenation
of $c$ and $c'$} is the $G(\Y)$--loop
$$c*c' =(g_0, e_1, \ldots, e_n, g_n g'_0, e'_1, \ldots, e'_{n'},
g'_{n'}).$$ The operation $[c][c']=[c*c']$ defines a group structure
on the set of homotopy classes of $G(\Y)$--loops at $\s_0$.

\begin{defn}[fundamental group of $G(\Y)$ at $\s_0$] The fundamental group of $G(\Y)$
at $\s_0$ is the set of homotopy classes of $G(\Y)$--loops at $\s_0$,
with the group structure induced by concatenation.  It is denoted
$\pi_1(G(\Y),\s_0)$.
\end{defn}

\noindent Different choices of basepoint $\s_0 \in V(\Y)$ result in
isomorphic fundamental groups.

The second definition of the fundamental group of a complex of
groups involves the choice of a maximal tree $T$ in the 1-skeleton
of the geometric realization $|\Y|$.  By abuse of notation, we will
say that $T$ is a maximal tree in $\Y$.

\begin{prop}[Theorem 3.7, \cite{BH}]  For any maximal tree $T$ in $\Y$, the fundamental group $\pi_1(G(\Y), \s_0)$ is isomorphic to the
abstract group $\pi_1(G(\Y), T)$, presented by the generating set
\[\coprod_{ \s \in V(\Y)} G_\s \coprod E^{\pm}(\Y)\] with the following relations:
\begin{enumerate}
 \item the relations in the groups $G_\s$;
\item $(a^+)^{-1}=a^-$ and $(a^-)^{-1}=a^+$;
\item  $a^+ b^+ = g_{a,b}(ab)^+$, for every $(a,b) \in E^{(2)}(\Y)$;
\item $\psi_a (g) = a^+ g a^-$, for every $g \in G_{i(a)}$; and
\item $a^+=1$ for every edge $a \in T$.\end{enumerate}
\end{prop}

\noindent If $\Y$ is simply connected, then $\pi_1(G(\Y),T)$ is
isomorphic to the direct limit of the diagram of groups $G_\s$ and
monomorphisms $\psi_a$.  The isomorphism $\pi_1(G(\Y), \s_0) \to
\pi_1(G(\Y), T)$ is the restriction of the natural projection
$FG(\Y) \to \pi_1(G(\Y), T)$.  Its inverse $\kappa_T$ is defined in
the proof of Proposition~\ref{p:isoms} in
Section~\ref{sss:univ_cover} below.

Let $\phi :G(\Y) \to G'(\Y')$ be a morphism over a morphism of
scwols $l: \Y \to \Y'$.  Then $\phi$ induces a homomorphism $F\phi:
FG(\Y) \to FG'(\Y')$, defined by $F\phi(g)=\phi_\s(g)$ for $g\in
G_\s$, and $F\phi(a^+)=\phi(a)l(a)^+$. The restriction of $F\phi$ to
$\pi_1(G(\Y), \s_0)$ is a natural homomorphism
\[\pi_1(\phi,\s_0):\pi_1(G(\Y), \s_0) \to \pi_1 (G'(\Y'), l(\s_0)).\]
In the particular case of a morphism $\phi:G(\Y) \to G$, where $G$
is a group, the induced homomorphism
\[\pi_1(\phi,\s_0):\pi_1(G(\Y),\s_0) \to G\] is defined by $g \mapsto
\phi_\s(g)$ for $g \in G_\s$, and $a^+ \mapsto \phi(a)$.

%*********************************************************************************
\subsubsection{Developments and the universal cover}\label{sss:univ_cover}
%*********************************************************************************

Any morphism from a complex of groups to a group induces a scwol,
called the development.

\begin{defn}[development]\label{d:development} Let $\phi: G(\Y) \to G$ be
a morphism from a complex of groups $G(\Y)$ to a group $G$.  The
scwol $D(\Y, \phi)$, called the development of $G(\Y)$ with respect
to $\phi$, is defined as follows.

The set of vertices is $$V(\D)=\{ ([g], \s): \s \in V(\Y), [g] \in
G/\phi_\s(G_\s) \}$$ and the set of edges is
$$E(\D)=\{ ([g], a) : a \in E(\Y), [g] \in G/\phi_{i(a)}(G_{i(a)}) \}.$$ The maps to initial and terminal vertices
are given by $$i([g], a)=([g], i(a))$$ and $$t([g],
a)=([g\phi(a)^{-1}], t(a))$$ and the composition of edges $([g],
a)([h],b)=([h], ab)$ is defined where $(a,b) \in E^{(2)}(\Y)$, $g, h
\in G$ and $g^{-1}h\phi(b)^{-1} \in \phi_{i(a)}(G_{i(a)})$.

 The
group $G$ acts naturally on $\D$: given $g, h \in G$ and $\alpha \in
\Y$, the action is $h \cdot ([g], \alpha)=([hg], \alpha)$.
\end{defn}

\begin{prop}[Theorems~2.13,~3.14 and~3.15,~\cite{BH}]\label{p:dev} Let $G(\Y)$
be a complex of groups over a connected scwol $\Y$ and let $G$ be a group.
\begin{enumerate} \item Let $\phi:G(\Y)\to G$ be a morphism which is injective
on the local groups. Then $G(\Y)$ is the complex of groups (with respect to
canonical choices) associated to the action of $G$ on the development $\D$, and
$\phi:G(\Y) \to G$ equals the canonical morphism $\phi_1:G(\Y) \to G$. \item
Suppose $G(\Y)$ is a complex of groups associated to the action of $G$ on a
simply connected scwol $\X$, and $\phi_1:G(\Y) \to G$ is the canonical
morphism.  Then $\phi_1$ induces a group isomorphism
\[\pi_1(\phi_1,\s_0):\pi_1(G(\Y),\s_0) \stackrel{\sim}{\to} G\] (see the final
paragraph of Section~\ref{sss:fundgp}), and there is a $G$-equivariant
isomorphism of scwols \[\Phi_1:D(\Y,\phi_1) \stackrel{\sim}{\to} \X\] given by,
for $g \in G$ and $\alpha \in \Y$, \[([g],\alpha) \mapsto g\cdot\overline\alpha.
\] \end{enumerate} \end{prop}

The following result, on the functoriality of developments, is used
to prove Theorem~\ref{t:functor_cov}, stated in the introduction.

\begin{prop}[Theorem 2.18,~\cite{BH}]\label{p:functor_dev} Let $G(\Y)$ and $G'(\Y')$ be complexes of groups over scwols
$\Y$ and $\Y'$.  Let $\phi:G(\Y) \to G$ and $\phi':G'(\Y') \to G'$
be morphisms to groups $G$ and $G'$ and let $\Lambda:G \to G'$ be a
group homomorphism. Let $\lambda:G(\Y) \to G'(\Y')$ be a morphism
over $l : \Y \to \Y'$.

Suppose
there is a homotopy from $\Lambda\phi$ to $\phi'\lambda$, given by elements $
k_\s \in G'$ (see Definition~\ref{d:homotopy}).  Then there is a
$\Lambda$-equivariant morphism of the developments
\[L:D(\Y,\phi) \to D(\Y',\phi')\]
given by, for $g \in G$ and $\alpha \in \Y$,
\[([g],\alpha) \mapsto ([\Lambda(g)k_{i(\alpha)}^{-1}], l(\alpha)).\]
Moreover, if $\phi$ and $\phi'$ are injective on the local groups,
and $\lambda$ and $\Lambda$ are isomorphisms, then $L$ is an
isomorphism of scwols.
\end{prop}

We now define the universal cover.

\begin{defn}[universal cover of a developable complex of
groups]\label{d:univ_cover} Let $G(\Y)$ be a developable complex of
groups over a connected scwol $\Y$. Choose a maximal tree $T$ in
$\Y$. Let
\[\iota_T: G(\mathcal{Y}) \to \pi_1(G(\mathcal{Y}),T)\] be the morphism of complexes of groups mapping the local group $G_\s$ to its image in $\fgYT$, and the edge $a$ to
the image of $a^+$ in $\fgYT$. The development
$\DYT=D(\Y, \iota_T)$ is called a universal cover of
 $G(\Y)$.
\end{defn}

\begin{thm}[Theorem~3.13,~\cite{BH}] The universal cover $\DYT$ is
connected and simply connected.
\end{thm}

\noindent As described in Definition~\ref{d:development}, the
fundamental group $\pi_1(G(\mathcal{Y}),T)$ acts canonically on
$\DYT$.

A group action on a scwol
induces the following explicit isomorphisms of groups and scwols.

\begin{prop}\label{p:isoms} Let $G$ be a group acting on a simply connected scwol $\mathcal{X}$, and let
$G(\Y)$ be the induced complex of groups (with respect to some
choices $C_\bullet=\{\overline\s,h_a\}$). Choose a maximal tree $T$
in $\Y$ and a vertex $\s_0 \in V(\Y)$.  For $e \in E^\pm(\Y)$, let
\[h_e = \left\{\begin{array}{ll} h_a & \mbox{ if }e = a^+
\\ h_a^{-1} & \mbox{ if }e = a^-\end{array}.\right. \]
For $\s \in V(\Y)$, let $c_\s = (e_1,e_2, \ldots,e_n)$ be the unique
edge-path contained in $T$, with no backtracking, which joins $\s_0$
to $\s$, and let $h_\s = h_{e_1}h_{e_2} \cdots h_{e_n}$.

Then there is a group isomorphism
\[\Lambda_T: \pi_1(G(\mathcal{Y}),T) {\to}  G \]
defined on generators by
\begin{align*}
g \mapsto & h_\sigma g h_{\sigma}^{-1} \mbox{ for } g \in G_\s\\
a^+ \mapsto & h_{t(a)} h_a h_{i(a)}^{-1}
\end{align*}
and a $\Lambda_T$-equivariant isomorphism of scwols
\begin{align*}
\tilde{L}_T: \DYT &\to \mathcal{X}\\
([g], \alpha) &\mapsto \Lambda_T(g) h_{i(\alpha)} \cdot
\overline{\alpha}.
\end{align*}
\end{prop}

\begin{pf}
For $\sigma \in V(\mathcal{Y})$ let $\pi_{\sigma}=e_1 e_2 \cdots e_n
$ be the element of $FG(\Y)$ corresponding to the edge-path $c_\s$.
Then by Theorem~3.7,~\cite{BH}, there is a canonical isomorphism
\[\kappa_T:\pi_1(G(\mathcal{Y}),T) \stackrel{\sim}{\to}
\pi_1(G(\mathcal{Y}), \sigma_0)\] defined on generators by
\begin{align*}
g \mapsto & \pi_\sigma g \pi_{\sigma}^{-1} \mbox{ for } g \in G_\s\\
a^+ \mapsto & \pi_{t(a)} a^+ \pi_{i(a)}^{-1}.
\end{align*}
By Proposition~\ref{p:dev}, the canonical morphism of complexes of
groups $\phi_1:G(\Y) \to G$ induces a group isomorphism
$\pi_1(\phi_1,\s_0):\pi_1(G(\Y),\s_0)\to G$.  Composing $\kappa_T$
with $\pi_1(\phi_1, \sigma_0)$, we obtain the group isomorphism
$\Lambda_T:\pi_1(G(\Y),T) \stackrel{\sim}{\to} G$ defined above.

We now have the square
$$\xymatrix{
G(\mathcal{Y}) \ar[r]^-{\iota_T} \ar[d]_{\lambda=Id} &
\pi_1(G(\mathcal{Y}),T) \ar[d]^-{\Lambda_T}\\
G(\mathcal{Y}) \ar[r]^-{\phi_1} & G. }
$$
This commutes up to a homotopy from $\Lambda_T \iota_T$ to
$\phi_1\lambda$, given by the elements $h_{\sigma}^{-1}$. Thus, by
Proposition~\ref{p:functor_dev}, there is a $\Lambda_T$-equivariant
morphism of scwols
\begin{align*}
L_T: D(\mathcal{Y},T) &\to D(\mathcal{Y}, \phi_1)\\
([g], \alpha) &\mapsto ([\Lambda_T(g) h_{i(\alpha)}], \alpha)
\end{align*}
which is an isomorphism since $\iota_T$ and $\phi_1$ are injective
on the local groups, and both $\lambda$ and $\Lambda_T$ are
isomorphisms. Composing $L_T$ with the $G$-equivariant isomorphism
$\Phi_1:D(\Y,\phi_1) \to \X$ (see Proposition~\ref{p:dev}), we
obtain a $\Lambda_T$-equivariant isomorphism of scwols
\begin{align*}
\tilde{L}_T: \DYT &\to \mathcal{X}\\
([g], \alpha) &\mapsto \Lambda_T(g) h_{i(\alpha)} \cdot
\overline{\alpha}
\end{align*} as required.\end{pf}

%*************************************************************************************************************
\subsubsection{Local developments and nonpositive curvature}\label{sss:nonpos}
%*************************************************************************************************************

Let $K$ be a connected polyhedral complex and let $\Y$ be the scwol
associated to $K$, so that $|\Y|$ is the first barycentric
subdivision of $K$. The \emph{star} $\St(\s)$ of a vertex $\s \in
V(\Y)$ is the union of the interiors of the simplices in $|\Y|$
which meet $\s$.  If $G(\Y)$ is a complex of groups over $\Y$, then
each $\s\in V(\Y)$ has a \emph{local development}, even if $G(\Y)$
is not developable. That is, we may naturally associate to each
vertex $\s \in V(\Y)$ an action of $G_\s$ on some simplicial complex
$\St(\tilde{\s})$ containing a vertex $\tilde\s$, such that
$\St(\s)$ is the quotient of $ \St(\tilde\s)$ by the action of
$G_\s$.  If $G(\Y)$ is developable, then for each $\s \in V(\Y)$,
the local development at $\s$ is isomorphic to the star of each lift
$\tilde\s$ of $\s$ in the universal cover $\DYT$.

We denote by $\st(\tilde\s)$ the star of $\tilde\s$ in
$\St(\tilde\s)$.

\begin{lem}[Lemma~5.2,~\cite{BH}]\label{l:star} Let $\lambda:G(\Y) \to G'(\Y')$ be a covering of
complexes of groups, over a morphism of scwols $l:\Y \to \Y'$.  Then
for each $\s \in V(\Y')$, Condition~\eqref{i:cov_bijection} in the
definition of a covering (Definition~\ref{d:covering}) is equivalent
to the existence of a $\lambda_\s$-equivariant bijection
$\st(\tilde\s) \to \st(\widetilde{l(\s)})$.
\end{lem}

In the case that $\Y$ is the scwol associated to a polyhedral
complex $K$, each local development $\St(\tilde\s)$ has a metric
structure induced by that of $K$ (see p. 562,~\cite{BH}).  A complex
of groups $G(\Y)$ has \emph{nonpositive curvature} if for all $\s
\in V(\Y)$, the local development at $\s$ has nonpositive curvature
(that is, $\St(\tilde\s)$ is locally CAT($\kappa$) for some $\kappa
\leq 0$) in this induced metric.  The importance of this condition
is given by:

\begin{thm}[Theorem~4.17,~\cite{BH}] If a complex of groups has nonpositive curvature, it is developable.\end{thm}

\noindent We will use the following condition to establish
nonpositive curvature:

\begin{lem}[Remark~4.18,~\cite{BH}]\label{l:CAT1} Let $\Y$ be the scwol associated to an $M_\kappa$-polyhedral complex $K$,
with $\kappa \leq 0$.  Then $G(\Y)$ has nonpositive curvature if and
only if, for each vertex $\tau$ of $K$, the geometric link of
$\tilde\tau$ in  $\st(\tilde\tau)$, with the induced spherical
structure, is CAT(1).\end{lem}

%*******************************************************************************************************
%*******************************************************************************************************
%*******************************************************************************************************
\section{Covering theory for complexes of groups}\label{s:covering_theory}
%*******************************************************************************************************
%*******************************************************************************************************
%*******************************************************************************************************

This section contains our results for complexes of groups which are analogous to those
for graphs of groups in~\cite{b1:ctgg}. We consider the functoriality of
morphisms of complexes of groups in Section~\ref{ss:functor_morphism} and that
of coverings in Section~\ref{ss:functor_cov}, culminating in the (constructive)
proof of  Theorem~\ref{t:functor_cov}. 
Section~\ref{ss:faithfulness} then characterizes faithfulness of complexes of
groups.  In Section~\ref{ss:main_lemma} a key technical result, the Main Lemma
(Lemma~\ref{l:main_lemma}),
is proved.  The Main Lemma makes precise the relationship between maps of groups and
scwols, and induced maps of fundamental groups and universal
covers of complexes of groups.  We consider the relationship between coverings
and developability in Section~\ref{ss:cov_dev}; this has no analogy for graphs
of groups since every graph of groups is developable.

%*******************************************************************************************************
%*******************************************************************************************************
\subsection{Functoriality of morphisms}\label{ss:functor_morphism}
%*******************************************************************************************************
%*******************************************************************************************************

Proposition~\ref{p:functor_morphism} below gives explicit
constructions of the maps on fundamental groups and universal covers
induced by a morphism of developable complexes of groups.

\begin{prop}\label{p:functor_morphism} Let $\lambda:G(\Y) \to G'(\Y')$
be a morphism of complexes of groups over a morphism of scwols $l:\Y
\to \Y'$, where $\Y$ and $\Y'$ are connected. Assume $G(\Y)$ and
$G'(\Y')$ are developable.  For any choice of $\s_0 \in V(\Y)$ and
maximal trees $T$ and $T'$ in $\Y$ and $\Y'$ respectively, $\lambda$
induces a homomorphism of fundamental groups
\[ \Lambda_{T,T'}= \Lambda^\lambda_{T,T'}:\fgYT \to \fgYpTp\]
and a $\Lambda_{T,T'}$-equivariant morphism of universal covers
\[L_{T,T'}^\lambda: \DYT \to \DYpTp.\]
\end{prop}

\begin{pf}  Let
$\s_0' = l(\s_0)$.  Recall from the proof of
Proposition~\ref{p:isoms} that there is a canonical isomorphism
\[\kappa_T: \pi_1(G(\Y),T) \stackrel{\sim}{\longrightarrow}
\pi_1(G(\Y),\sigma_0)\] and from the last paragraph of
Section~\ref{sss:fundgp} that the morphism $\lambda$ induces a group
homomorphism $\pi_1(\lambda, \sigma_0):\pi_1(G(\Y),\s_0) \to
\pi_1(G'(\Y'),\s_0')$ which is the restriction of the morphism
$F\lambda :FG(\Y) \to FG'(\Y')$. The group homomorphism
\[\Lambda_{T,T'}: \pi_1(G(\Y), T) \to \pi_1(G'(\Y'),T')\] is defined by
the composition $\kappa'^{-1}_{T'} \circ \pi_1(\lambda,\s_0) \circ
\kappa_T$:
$$ \pi_1(G(\Y), T) \stackrel{\sim}{\longrightarrow} \pi_1(G(\Y), \sigma_0) \longrightarrow \pi_1(G'(\Y'),
\sigma_0') \stackrel{\sim}{\longrightarrow} \pi_1(G'(\Y'), T')$$

We now have a square
$$\xymatrix{
G(\Y) \ar[r]^-{\iota_T} \ar[d]_-{\lambda} & \pi_1(G(\Y), T)
\ar[d]_-{\Lambda_{T,T'}} \\
G'(\Y') \ar[r]^-{\iota'_{T'}} & \pi_1(G'(\Y'), T').}$$ We claim that
that there is a homotopy from $\Lambda_{T,T'} \circ \iota_{T}$ to
$\iota'_{T'} \circ \lambda$. For $\sigma \in V(\Y)$ let $\pi_\s=e_1
e_2 \cdots e_n$ be the element of $FG(\Y)$ corresponding to the
unique path $(e_1,e_2,\ldots,e_n)$ in $T$ without backtracking from
$\s_0$ to $\s$, and similarly for $\pi'_{l(\s)} \in FG'(\Y')$. Then
for $g \in G_\s$, we have
\begin{equation*}
\begin{split}
(\Lambda_{T,T'} \circ \iota_T) (g)=&\Lambda_{T,T'}(g)\\
=&\kappa'^{-1}_{T'} \circ \pi_1(\lambda,\sigma_0) \circ
\kappa_T(g)\\
=&\kappa'^{-1}_{T'} \circ \pi_1(\lambda, \sigma_0)(\pi_{\sigma}
g \pi_{\sigma}^{-1})\\
=&\kappa'^{-1}_{T'}\{F\lambda(\pi_\s) \lambda_{\sigma}(g) (F\lambda(\pi_\s))^{-1}\}\\
=&\kappa'^{-1}_{T'}\{ F\lambda(\pi_\s)( \pi'_{l(\sigma)})^{-1}\}\,
(\iota'_{T'}\! \circ\lambda_{\sigma})(g)\,
\kappa'^{-1}_{T'}\{\pi'_{l(\sigma)} (F\lambda(\pi_\s))^{-1} \}.
\end{split}
\end{equation*}
Setting
\[u_{\sigma}=\kappa'^{-1}_{T'}\{F\lambda(\pi_\s) (\pi'_{l(\sigma)})^{-1}\} \in \fgYpTp\] we
conclude
$$(\Lambda_{T,T'} \circ \iota_T)(g)=u_{\sigma}
(\iota'_{T'} \circ
\lambda_{\sigma})(g)\,u_{\sigma}^{-1}=\Ad(u_{\sigma})(\iota'_{T'}
\circ \lambda)(g).$$ Similarly, if $a \in E(\Y)$, we compute
\begin{equation*}
\begin{split}
(\Lambda_{T,T'} \circ \iota_{T})(a)&=\kappa'^{-1}_{T'} \circ
\pi_1(\lambda, \sigma_0) \circ \kappa_{T}(a^+)
\\&=\kappa'^{-1}_{T'} \circ \pi_1(\lambda, \sigma_0)(\pi_{t(a)} a^+
\pi_{i(a)}^{-1})\\&=u_{t(a)} \lambda(a) l(a)^+ u_{i(a)}^{-1}
\\&=u_{t(a)} (\iota'_{T'} \circ \lambda)(a) u^{-1}_{i(a)}.
\end{split}
\end{equation*} The last equality comes from the definition of composition of morphisms,
\[(\iota'_{T'} \circ
\lambda)(a)=(\iota'_{T'})_{l(t(a))}(\lambda(a))\iota'_{T'}(l(a))=\lambda(a)
l(a)^{+}.\]
  Hence the desired homotopy from $\Lambda_{T,T'} \circ \iota_{T}$ to
$\iota'_{T'} \circ \lambda$ is given by the elements
$u^{-1}_{\sigma}$.

By Proposition~\ref{p:functor_dev} there is thus a
$\Lambda_{T,T'}$-equivariant morphism of universal covers
$$L^\lambda_{T,T'}: D(\Y, T) \to D(\Y', T')$$
given by
$$([g], \alpha) \mapsto ([\Lambda_{T,T'}(g)
u_{i(\alpha)}], l(\alpha)).$$
\end{pf}

Corollary~\ref{c:commute} below says that if a diagram of morphisms
of developable complexes of groups commutes, then the corresponding
diagrams of the induced maps on fundamental groups and universal
covers, defined in Proposition~\ref{p:functor_morphism} above, also
commute.

\begin{cor}\label{c:commute}  With the notation of Proposition~\ref{p:functor_morphism}, let $G''(\Y'')$ be a
developable complex of groups over a connected scwol $\Y''$, and
assume there is a morphism $\lambda':G'(\Y') \to G''(\Y'')$.  Choose
a maximal tree $T''$ in $\Y''$. Then the composition
\[\lambda''=\lambda' \circ \lambda\] induces a group homomorphism
$\Lambda_{T,T''}: \pi_1(G(\Y), T) \to \pi_1(G''(\Y''), T'')$ and a
$\Lambda_{T,T''}$-equivariant morphism of universal covers
$L^{\lambda''}_{T,T''}: D(\Y, T) \to D(\Y'', T'')$, such that
\[ L^{\lambda''}_{T,T''}=L^{\lambda'}_{T',T''} \circ
L^\lambda_{T,T'}\] and \[\Lambda_{T,T''}=\Lambda_{T',T''} \circ
\Lambda_{T, T'}.\]
\end{cor}

\begin{pf} The proof follows from the constructions given in
Proposition~\ref{p:functor_morphism} above, and the definition of
composition of morphisms.
\end{pf}

%*******************************************************************************************************
%*******************************************************************************************************
\subsection{Functoriality of coverings}\label{ss:functor_cov}
%*******************************************************************************************************
%*******************************************************************************************************

In this section we prove Theorem~\ref{t:functor_cov}, stated in the
Introduction.  The maps $\Lambda_{T,T'}$ and $L^\lambda_{T,T'}$ are
those defined in Proposition~\ref{p:functor_morphism} above.

\begin{prop}\label{p:induced} Let $\lambda:G(\cY)\to G'(\cY')$ be a covering of complexes of groups over a morphism of scwols $l:\Y\to\Y'$, where
$\Y$ and $\Y'$ are connected.  Assume $G(\Y)$ and $G'(\Y')$ are
developable.  For any choice of $\s_0 \in V(\Y)$ and maximal trees
$T$ and $T'$ in $\Y$ and $\Y'$ respectively, the induced
homomorphism of fundamental groups
\[\Lambda_{T,T'}: \fgYT \rightarrow \fgYpTp\] is a monomorphism and
\[L^\lambda_{T,T'}:D(\Y,T) \rightarrow D(\Y',T') \] is a
$\Lambda_{T,T'}$-equivariant isomorphism of scwols.
\end{prop}

\begin{pf} We begin with Lemma~\ref{l:covmorphism} below, which shows that $L^\lambda_{T,T'}$ is
a covering of scwols (see Definition~\ref{d:cov_scwols}).
Corollary~\ref{c:isomscwols} of this lemma shows that
$L^\lambda_{T,T'}$ is an isomorphism of scwols. We then use this
result to show that $\Lambda_{T,T'}$ is injective.

\begin{lem}\label{l:covmorphism} The morphism $L^\lambda_{T,T'}$ is a covering of scwols.  \end{lem}

\begin{pf}
 Let $g \in \fgYT$ and $\sigma \in V(\cY)$.

 We first show that $L^\lambda_{T,T'}$ is injective on the set of
edges with terminal vertex $([g],\sigma)$. Suppose $a_1$ and $a_2$
are edges of $\cY$ (with $t(a_1) = t(a_2) = \sigma$), that for some
$h_1,h_2 \in \fgYT$  \[ t\left([h_1],a_1\right) =
([g],\sigma)=t\left([h_2],a_2\right) \] and that
\[L^\lambda_{T,T'}\left([h_1],a_1\right)=
L^\lambda_{T,T'}\left([h_2],a_2\right).\] By definition of
$L^\lambda_{T,T'}$, we then have $l(a_1)=l(a_2) = a'$ say, with
$t(a') = l(t(a_1))= l(\sigma)=\sigma'$.  Also, by definition of the
map $t:E(\DYT)\to V(\DYT)$, we have, for some $h \in G_\s$, \[ h_1
{a_1^-} = h_2 {a_2^-} h^{-1}.
\] Now by definition of $L_{T,T'}^\lambda$, it follows that the group $G'_{i(a')}$ contains
\begin{equation*}
\begin{split} & \left(\Lambda_{T,T'}  (h_1)u_{i(a_1)} \right)^{-1}
 \left(\Lambda_{T,T'}(h_2)u_{i(a_2)} \right)
\\&=  u_{i(a_1)}^{-1} \Lambda_{T,T'}\left(a_1^-\,h\,a_2^+\right) u_{i(a_2)} \\
&=  u_{i(a_1)}^{-1}  u_{i(a_1)} l(a_1)^- \lambda(a_1)^{-1}
u_{t(a_1)}^{-1} u_\s \lambda_\s(h) u_\s^{-1} u_{t(a_2)} \lambda(a_2)
l(a_2)^+ u_{i(a_2)}^{-1} u_{i(a_2)} \\  &=a'^- \lambda(a_1)^{-1}
\lambda_\s(h) \lambda(a_2) a'^+.\end{split}
\end{equation*} Thus by the relation $a'^+ k a'^-=\psi_{a'}(k)$, for all $k \in G'_{i(a')}$,
\[ \lambda(a_1)^{-1}\, \lambda_\s(h)\, \lambda(a_2) \in
\psi_{a'}(G'_{i(a')}). \] That is, $\lambda(a_1)$ and $\lambda
_\sigma(h) \lambda(a_2)$ belong to the same coset of
$\psi_{a'}(G'_{i(a')})$ in $G'_{\sigma'}$. By
Condition~\eqref{i:cov_bijection} in the definition of a covering
(Definition~\ref{d:covering}, this implies $a_1 = a_2 = a$, say, and
$h \in \psi_a(G_{i(a)})$.  It follows that $h_1$ and $h_2$ belong to
the same coset of $G_{i(a)}$ in $\fgYT$. Thus $L^\lambda_{T,T'}$ is
injective on the set of edges with terminal vertex $([g],\sigma)$.

We now show that $L^\lambda_{T,T'}$ surjects onto the set of edges
of $\DYpTp$ with terminal vertex $L^\lambda_{T,T'}([g],\sigma)$.
Suppose  \[t\left([h'],a'\right) =
L^\lambda_{T,T'}\left([g],\sigma\right)\] where $h' \in \fgYpTp$,
$a' \in E(\cY')$.  Then $t(a') = \sigma'=l(\s)$ and by definition of
$L_{T,T'}^\lambda$,
\begin{equation}\label{e:hp} h' a'^- = \Lambda_{T,T'}(g) u_\s
k_{\s'} \end{equation} for some $k_{\sigma'} \in G'_{\sigma'}$. By
Condition~\ref{i:cov_bijection} in the definition of a covering,
there exists an edge $a \in E(\cY)$ with $l(a)=a'$ and $t(a) =
\sigma$, and an element $k_\sigma \in G_\sigma$, such that $ \lambda
_\sigma(k_\sigma) \lambda (a)$ and $k_{\sigma'}$ belong to the same
coset of $\psi_{a'}(G_{i(a')})$ in $G'_{\sigma'}$. Let $h = g
k_\sigma a^+ \in \fgYT$ and note that by
Definition~\ref{d:development},
\[
t\left([h],a\right) = \left([gk_\s a^+ \iota_T(a)^{-1}],
t(a)\right)= \left([gk_\sigma a^+ a^-],
\sigma\right)=\left([gk_\sigma], \sigma\right) =
\left([g],\sigma\right).
\]
We claim
\[
L^\lambda_{T,T'}\left([h],a\right) = \left([h'],a'\right).
\]
By Equation~\eqref{e:hp} above, the choice of $a$ and $k_\sigma$ and
the relation $\psi_{a'}(k')=a'^+ k' a'^-$ for all $k' \in
G'_{i(a')}$, we have
\begin{align*}
\Lambda_{T,T'}(h) u_{i(a)} &=
 \Lambda_{T,T'}(g)u_\s \lambda_\s(k_\s) u_\s^{-1} u_{t(a)} \lambda(a)
 l(a)^+ u_{i(a)}^{-1} u_{i(a)} \\ &=
 h' a'^- k_{\sigma'}^{-1}
 \lambda _\sigma(k_\sigma)
 \lambda(a) a'^+ \\ &\in h'G'_{i(a')}.
\end{align*}
Hence,
\[L^\lambda_{T,T'}\left( [h], a \right)=
\left( \left[\Lambda_{T,T'}(h) u_{i(a)}\right], a'\right) =
([h'],a').\] We conclude that $L^\lambda_{T,T'}$ is a covering of
scwols.\end{pf}

\begin{cor}\label{c:isomscwols} Under the assumptions of Proposition~\ref{p:induced},
the morphism $L^\lambda_{T,T'}:\DYT\to\DYpTp$ is an isomorphism of scwols.
\end{cor}

\begin{pf} By Lemma~\ref{l:covmorphism}, $L^\lambda_{T,T'}$ is a covering
morphism.  Since $\DYpTp$ is connected, $L^\lambda_{T,T'}$ is
surjective, and since $\DYT$ is connected and $\DYpTp$ is simply
connected, $L^\lambda_{T,T'}$ is injective. See Remark
1.9(2),~\cite{BH}. \end{pf}

We complete the proof of Proposition~\ref{p:induced} by showing that
$\Lambda_{T,T'}$ is a monomorphism of groups. Suppose $g \in \fgYT$
and $\Lambda_{T,T'}(g)=1$.  Since $L^\lambda_{T,T'}$ is injective
and $\Lambda_{T,T'}$-equivariant, $g$ must act trivially on $\DYT$.
In particular,
\[
g \cdot \left( [1], \sigma_0 \right) = \left( [g], \sigma_0 \right)
= \left( [1], \sigma_0 \right)
\]
so $g \in G_{\s_0}$.  We then calculate
\begin{align*}
\Lambda_{T,T'}(g) &= \kappa'^{-1}_{T'} \circ \pi_1(\lambda,\s_0) \circ \kappa_T ((\iota_T)_{\s_0}(g)) \\
&= \kappa'^{-1}_{T'} (\lambda_{\s_0}((\iota_T)_{\s_0}(g)))\\
&= 1.
\end{align*}
Since $\kappa'^{-1}_{T'}$, $\lambda _{\sigma_0}$ and
$(\iota_T)_{\s_0}$ are each injective, this implies $g=1$. Thus
$\Lambda_{T,T'}$ is injective.\end{pf}

\begin{cor}\label{c:induced} Let $\lambda:G(\cY)\to G'(\cY')$ be a covering of complexes of groups.  Suppose for some $\kappa \in \R$ that the scwols $\cY$ and $\cY'$ are
associated to $M_\kappa$-polyhedral complexes with finitely many
isometry classes of cells. If $G(\cY)$ and $G'(\cY')$ are
developable, then the geometric realizations of their respective
universal covers are isometric (as polyhedral complexes).
\end{cor}

%*******************************************************************************************************
%*******************************************************************************************************
\subsection{Faithfulness}\label{ss:faithfulness}
%*******************************************************************************************************
%*******************************************************************************************************

\begin{defn}[faithful]\label{d:faithful} Let $G(\Y)$ be a developable complex of groups.  We say $G(\Y)$ is faithful if the natural homomorphism
$\fgYT \to \Aut(\DYT)$ is a monomorphism, for any choice of maximal
tree $T$ in $\Y$.
\end{defn}

If $G(\Y)$ is a complex of groups associated to the action of a
group $G$ on a scwol $\X$, then $G(\Y)$ is faithful.

Proposition~\ref{p:NT} below may be used to give sufficient
conditions for faithfulness.

\begin{prop}\label{p:NT} Let $G(\cY)$ be a
developable complex of groups over a connected scwol $\cY$. Choose a
maximal tree $T$ in $\cY$, and identify each local group $G_\s$ with
its image in $\fgYT$ under the morphism $\iota_T$.  Let
\[N_T = \ker(\pi_1(G(\cY),T) \rightarrow \DYT).\] Then
\begin{enumerate} \item \label{i:vertex} $N_T$ is a vertex
subgroup, that is $N_T \leq {G}_\sigma$ for each $\sigma \in
V(\cY)$. \item \label{i:inv} $N_T$ is $\cY$-invariant, that is
$\psi_a(N_T) = N_T$ for each $a \in E(\cY)$. \item \label{i:normal}
$N_T$ is normal, that is $N_T \unlhd {G}_\sigma$ for each $\sigma
\in V(\cY)$.\item\label{i:max} $N_T$ is maximal: if $N'_T$ is
another $\Y$-invariant normal vertex subgroup then $N_T' \leq N_T$.
\end{enumerate} \end{prop}

\begin{pf} If $h \in N_T$, then for all $\sigma \in V(\cY)$,
\[h \cdot ([1],\sigma) =
([h],\sigma) = ([1],\sigma)\] thus $h \in {G}_\sigma$. This
proves~\eqref{i:vertex}. Since $N_T$ is normal in $\pi_1(G(\cY),T)$
it is normal in each ${G}_\sigma$, proving~\eqref{i:normal}.

To prove~\eqref{i:inv}, let $a \in E(\Y)$.  In the group
$\pi_1(G(\cY), T)$ the following relation holds for each $g \in
{G}_{i(a)}$: \[\psi_a(g) = a^+ g a^-.\] Since $N_T$ is a subgroup of
${G}_{i(a)}$ and $N_T$ is normal in $\pi_1(G(\cY),T)$, it follows
that \[\psi_a(N_T) = a^+ N_T a^- = N_T\] as required.

To prove~\eqref{i:max}, we have, for all $g \in \fgYT$ and $\alpha
\in \Y$,
\[N_T' \cdot ([g],\alpha) = gN_T'g^{-1} \cdot ([g],\alpha) = g \cdot ([1],\alpha) = ([g],\alpha)\]
since $N_T'$ is normal in $\fgYT$ and $N'_T$ is a subgroup of
$G_{i(\alpha)}$.  Hence $N_T'$ is contained in $N_T$, as claimed.
\end{pf}

%******************************************************************************************************
%******************************************************************************************************
\subsection{Other functoriality results}\label{ss:main_lemma}
%******************************************************************************************************
%******************************************************************************************************

This section contains results similar to those in
Section~4,~\cite{b1:ctgg}.

We first prove the following useful characterization of isomorphisms
of complexes of groups. This result corresponds to
Corollary~4.6,~\cite{b1:ctgg}.

\begin{prop}\label{p:functor_isom}  Let $\lambda:G(\Y) \to G'(\Y')$ be a morphism of
developable complexes of groups over a morphism of scwols $l:\Y \to
\Y'$, where $\Y$ and $\Y'$ are connected scwols.  For any choice of
$\s_0 \in V(\Y)$ and maximal trees $T$ and $T'$ in $\Y$ and $\Y'$
respectively, $\lambda$ is an isomorphism if and only if both of the
maps $L_{T,T'}^\lambda$ and $\Lambda_{T,T'}$ are isomorphisms.
\end{prop}

\begin{pf}  If $\lambda$ is an isomorphism, it is clearly a
covering.  Proposition~\ref{p:induced} thus implies that
$L_{T,T'}^\lambda$ is an isomorphism of scwols and $\Lambda_{T,T'}$
is a monomorphism of groups.  Since $\lambda^{-1}$ is also a
covering, $\Lambda_{T,T'}^{-1} = (\Lambda_{T,T'})^{-1}$ is also a
monomorphism, hence $\Lambda_{T,T'}$ is an isomorphism.

Conversely, suppose $\lambda$ is not an isomorphism, thus one of
$\lambda$ and $\lambda^{-1}$ is not a covering.  Without loss of
generality, we assume $\lambda$ is not a covering.  Then either
\begin{enumerate}
\item\label{i:local} there is a homomorphism $\lambda_\s : G_\s \to G'_{l(\s)}$ which
is not injective, or
\item\label{i:star} there exists $a' \in
E(\Y')$ and $\s \in V(\Y)$ with $t(a')=\s'=l(\s)$, such that the map
\[
\coprod_{\substack{a \in l^{-1}(a') \\ t(a)=\s}} G_\s
/\psi_{a}(G_{i(a)}) \to G'_{\s'}/\psi_{a'} (G'_{i(a')})\] induced by
\[
 g \mapsto \lambda_{\s}(g) \lambda(a)
 \] is not bijective.
\end{enumerate}

Condition~\eqref{i:local} implies that the map $\Lambda_{T,T'}$ is
not a monomorphism at $G_\s$, thus $\Lambda_{T,T'}$ is not an
isomorphism.  Condition~\eqref{i:star} implies that
$L_{T,T'}^\lambda$ is not a local bijection at $\St(\tilde\s)$ (see
Remark~5.3,~\cite{BH}), thus the map $L_{T,T'}^\lambda$ is not an
isomorphism.
\end{pf}

The Main Lemma below, which corresponds to
Proposition~4.4,~\cite{b1:ctgg}, will be used many times in
Section~\ref{s:bijection}.  The data for the Main Lemma is as
follows.

Let $\X$ and $\X'$ be simply connected scwols, acted upon by groups
$G$ and $G'$ respectively, with quotient scwols $\Y = G\bs\X$ and
$\Y'=G' \bs \X'$.  Let $G(\Y)_{C_\bullet}$ and
$G'(\Y')_{C'_\bullet}$ be complexes of groups associated to the
actions of $G$ and $G'$, with respect to choices
$C_\bullet=(\overline\s,h_a)$ and
$C'_\bullet=(\overline{\s'},h_{a'})$.

Suppose $L: \mathcal{X} \to \mathcal{X}'$ is a morphism of scwols
which is equivariant with respect to some group homomorphism
$\Lambda: G \to G'$.  Let $l:\Y \to\Y'$ be the induced morphism of
quotient scwols.  Fix $\sigma_0 \in \mathcal{Y}$ and let
$\s_0'=l(\s_0)$.   Let $N_\bullet=\{k_\s\}$ be a set of elements of
$G'$ such that $k_\s\cdot L(\overline\s) = \overline{l(\s)}$ for all
$\s \in V(\Y)$.

With respect to these choices, there is an induced morphism
$\lambda=\lambda_{C_\bullet,C'_\bullet,N_\bullet}: G(\Y) \to
G'(\Y')$ (see Definition~\ref{d:induced_morphism}). For any choice
of maximal trees $T$ and $T'$ in $\Y$ and $\Y'$, respectively, let
$$\Lambda^\lambda_{T,T'}: \fgYT \to \fgYpTp$$
be the homomorphism of groups induced by $\lambda$ and let
$$L^\lambda_{T,T'}: \DYT \to \DYpTp$$
be the associated $\Lambda^\lambda_{T,T'}$-equivariant morphism of
scwols (see Proposition~\ref{p:functor_morphism}).  By
Proposition~\ref{p:isoms} we have isomorphisms of scwols
\[\tilde{L}_T: \DYT \stackrel{\sim}{\longrightarrow} \mathcal{X}\quad
\mbox{and}\quad \tilde{L}_{T'}: \DYpTp
\stackrel{\sim}{\longrightarrow} \mathcal{X}'\] which are
equivariant with respect to group isomorphisms
\[\Lambda_T: \fgYT \stackrel{\sim}{\longrightarrow} G\quad
\mbox{and}\quad \Lambda_{T'}: \fgYpTp
\stackrel{\sim}{\longrightarrow} G'\] respectively.

\begin{lem}[Main Lemma]\label{l:main_lemma} Suppose $C_\bullet$ and $C_\bullet'$ are
chosen so that $L(\overline{\s_0})=
\overline{l(\s_0)}=\overline{\s_0'}$, and $N_\bullet$ is chosen so
that $k_{\s_0}=1$.  Then the following diagrams commute:
\begin{enumerate} \item\label{i:commute_group}
$$\xymatrix{ \fgYT \ar[rr]^-{\Lambda^\lambda_{T,T'}}
\ar[d]^-{\Lambda_T} & & \fgYpTp
\ar[d]^-{\Lambda_{T'}}\\
G \ar[rr]^-{\Lambda} & &G'}$$
\item\label{i:commute_scwol}
$$\xymatrix{
\DYT \ar[rr]^-{L^\lambda_{T,T'}}  \ar[d]^-{\tilde{L}_T} & &\DYpTp
\ar[d]^-{\tilde{L}_{T'}}\\
\mathcal{X} \ar[rr]^-{L} && \mathcal{X}'.}$$\end{enumerate}
\end{lem}

\begin{pf} We first show the commutativity of~\eqref{i:commute_group}, and then use this
diagram and equivariance to prove that~\eqref{i:commute_scwol}
commutes.

By construction, \[\Lambda_T= \pi_1(\phi_1, \sigma_0) \circ
\kappa_T\quad \mbox{and}\quad \Lambda_{T'}= \pi_1(\phi'_1,
\sigma'_0) \circ \kappa'_{T'}\] where $\phi_1:G(\Y) \to G$ and
$\phi'_1:G'(\Y') \to G'$ are the canonical morphisms. Also,
$\Lambda^\lambda_{T,T'} = \kappa'^{-1}_{T'} \circ \pi_1 (\lambda,
\sigma_0) \circ \kappa_T$. Therefore it is enough to show that the
following diagram commutes:
$$\xymatrix{ \pi_1(G(\Y),\s_0) \ar[rr]^-{\pi_1(\lambda,\s_0)}
\ar[d]^-{\pi_1(\phi_1,\s_0)} & &\pi_1(G'(\Y'),\s_0')
\ar[d]^-{\pi_1(\phi_1',\s_0')}\\
G \ar[rr]^-{\Lambda} & &G'.}$$ Let $x \in \pi_1(G(\Y),\s_0)$.  Then
$x$ has the form
\[x=g_{\sigma_0}e_1g_{\s_1} \cdots e_n g_{\sigma_n}\] where
$(g_{\s_0},e_1,g_{\s_1},\ldots,e_n,g_{\s_n})$ is a $G(\Y)$--loop
based at $\s_0=\s_n$.  It follows that
\[\pi_1(\phi_1,\s_0)(x) = g_{\s_0} h_{e_1} g_{\s_1}  \cdots h_{e_n} g_{\s_n}\]
where the elements $h_{e_j}$ are as defined in
Proposition~\ref{p:isoms}. We now compute
\begin{equation*}
\begin{split}
\pi_1&(\phi_1', \sigma_0') \circ \pi_1(\lambda, \sigma_0)(x)\\
&=(k_{\sigma_0}\Lambda(g_{\sigma_0})k_{\sigma_0}^{-1})(k_{\sigma_0}\Lambda(h_{e_1})k_{\sigma_1}^{-1}h^{-1}_{l(e_1)})
h_{l(e_1)} (k_{\sigma_1} \Lambda(g_{\sigma_1}) k_{\sigma_1}^{-1})
\cdots
(k_{\sigma_n}\Lambda(g_{\sigma_n}) k_{\sigma_n}^{-1})\\
&=k_{\sigma_0}\Lambda(g_{\sigma_0}h_{e_1}g_{\sigma_1} \cdots
h_{e_n}g_{\sigma_n})
k_{\sigma_n}^{-1}\\
&=\Lambda\circ\pi_1(\phi_1,\s_0)(x)
\end{split}
\end{equation*}
since $k_{\sigma_0}=k_{\sigma_n}=1$. Thus~\eqref{i:commute_group}
commutes.

To prove that~\eqref{i:commute_scwol} commutes, let
\[\tilde{L} = \tilde{L}_{T'} \circ L^\lambda_{T,T'} \circ \tilde{L}_{T}^{-1}\]
We will show that $\tilde{L} = L$.  By the equivariance of the
morphisms of scwols used to define $\tilde{L}$, and the
commutativity of~\eqref{i:commute_group}, we have that $\tilde{L}$
is $\Lambda$-equivariant.  Thus it is enough to check (for example)
that $\tilde{L}(h_{i(\alpha)}\overline{\alpha})=
{L}(h_{i(\alpha)}\overline{\alpha})$ for all $\alpha \in \Y$. By
Proposition~\ref{p:isoms},
\begin{equation*}
\begin{split}
\tilde{L}(h_{i(\alpha)}\overline{\a})&=\tilde{L}_{T'} \circ L^\lambda_{T,T'}([1],\a)\\
&=\tilde{L}_{T'}([u_{i(\a)}],l(\a))\\
&=\Lambda_{T'}(u_{i(\a)})h_{i(l(\a))}\overline{l(\a)}.
\end{split}
\end{equation*}
 Let $\pi_{i(\alpha)}=e_1 e_2\cdots e_n$
be the element of $FG(\Y)$ which corresponds to the non-backtracking
path in $T$ from $\s_0$ to $i(\alpha)$, and similarly for $\pi
'_{i(l(\alpha))}=e_1'e_2'\cdots e'_{n'}$ in $FG'(\Y')$.  Then
\begin{equation*}\label{E:vroum}
\begin{split}
\Lambda_{T'}(u_{i(\a)})&=\Lambda_{T'} \circ
\kappa'^{-1}_{T'}\left\{F\lambda(\pi_{i(\a)})(\pi'_{i(l(\a))})^{-1}\right\}\\
&=\pi_1(\phi'_1,
\sigma_0')\left\{F\lambda(\pi_{i(\a)})(\pi'_{i(l(\a))})^{-1}\right\}\\
&=\pi_1(\phi'_1, \sigma_0')\left\{F\lambda(e_1)F\lambda(e_2)\cdots
F\lambda(e_n)e'^{-1}_{n'}\cdots e'^{-1}_2 e'^{-1}_1\right\} \\ &=
k_{\s_0} \Lambda(h_{e_1}) k_{\s_1}^{-1} k_{\s_1} \Lambda(h_{e_2})
k_{\s_2}^{-1} \cdots k_{\s_{n-1}} \Lambda(h_{e_n})
k_{\s_n}^{-1}h^{-1}_{e'_{n'}}\cdots h^{-1}_{e'_{2}} h^{-1}_{e'_{1}}
\\ &= k_{\s_0} \Lambda(h_{e_1} h_{e_2} \cdots h_{e_n}) k_{\s_n}^{-1}
(h_{e_1'} h_{e'_2} \cdots h_{e'_{n'}})^{-1} \\ &=
\Lambda(h_{i(\alpha)}) k_{\s_n}^{-1} h_{i(l(\a))}^{-1}
\end{split}
\end{equation*} since $k_{\sigma_0}=1$. Substituting, we obtain finally
\begin{equation*}
\begin{split}
\tilde{L}(h_{i(\a)}\overline{\a})&=\Lambda(h_{i(\a)})
k^{-1}_{\s_n}\overline{l(\a)}\\
&=\Lambda(h_{i(\a)}) k^{-1}_{i(\a)} \overline{l(\a)}\\
&=\Lambda(h_{i(\a)}) L(\overline{\a})\\
&=L(h_{i(\a)}\overline{\a})
\end{split}
\end{equation*}
as desired. This completes the proof of the Main Lemma.\end{pf}

The following result makes precise the relationship between a
developable complex of groups $G(\Y)$ and the complex of groups
induced by the action of $\fgYT$ on $\DYT$, for some maximal tree
$T$ in $\Y$.  It will be used to prove the Corollary to the Main Lemma below.

\begin{lem}\label{l:useful}  Let $G(\Y)$ be a developable complex of groups over a connected
scwol $\Y$. Choose a vertex $\s_0 \in V(\Y)$ and a maximal tree $T$
in $\Y$. Let $\Z$ be the quotient scwol
\[ \Z = \fgYT \bs \DYT \]
and let $f$ be the canonical isomorphism of scwols \begin{align*}
f: \Y &\to \Z \\
\alpha &\mapsto \fgYT \cdot ([1],\alpha)
\end{align*}
 Let $C_\bullet$ be the following
data for the action of $\fgYT$ on $\DYT$:
\[ \overline{f(\alpha)} = ([1],\alpha) \quad\mbox{and}\quad h_{f(a)} =
a^+\] and let $G(\Z)_{C_\bullet}$ be the complex of groups
associated to this data.  Then there is an isomorphism of complexes
of groups
\[\theta:G(\Y) \to G(\Z)\]
over $f$ such that
\[ \Lambda^\theta_{T,f(T)} = \Lambda_{f(T)}^{-1} \quad\mbox{and}\quad
L^\theta_{T,f(T)} = \tilde{L}_{f(T)}^{-1} \]
where $f(T)$ is the image of $T$ in $\Z$.
\end{lem}

\begin{pf}  We define $\theta$ by $\theta_\s(g) =
g$ for each $g \in G_\s$, and $\theta(a) = 1$ for each $a \in E(\Y)$
(here we are identifying $G_\s$ with its image in $\fgYT$).

 We then
have
\[ \Lambda^\theta_{T,f(T)} \circ \Lambda_{f(T)} = \kappa_{f(T)}^{-1} \circ \pi_1(\theta,\s_0) \circ \kappa_T \circ \pi_1(\phi_1,f(\s_0)) \circ
\kappa_{f(T)}. \] We claim that
\begin{equation}\label{e:useful} \pi_1(\theta,\s_0) \circ \kappa_T
\circ \pi_1(\phi_1,f(\s_0)) = 1.\end{equation} Let $g  \in
\pi_1(G(\Z),f(\s_0))$. Then $g = g_0 f(e_1) g_1 \cdots f(e_n) g_n$
for some $G(\Z)$-loop $(g_0,f(e_1),g_1,\ldots,f(e_n),g_n)$ based at
$f(\s_0)=f(\s_n)$,
and so
\begin{align*}
\kappa_T\circ \pi_1(\phi_1,f(\s_0))(g) &= \kappa_T( g_0 h_{f(e_1)}
g_1 \cdots h_{f(e_n)} g_n) \\ &= \pi_{\s_0} g_0 \pi_{\s_0}^{-1}
\kappa_T(h_{f(e_1)}) \pi_{\s_1} g_1 \pi_{\s_1}^{-1} \cdots
\kappa_T(h_{f(e_n)}) \pi_{\s_n} g_n \pi_{\s_n}^{-1}
\end{align*}
where $\pi_\s$ is the unique non-backtracking path in $T$ from
$\s_0$ to $\s$.  Now, applying $h_{f(a)} = a^+$
and $\kappa_T(a^+) = \pi_{t(a)}a^+ \pi_{i(a)}^{-1}$, as well as
$\pi_{\s_0} = \pi_{\s_n} = 1$, we have
\begin{align*}
\pi_1(\theta,\s_0) \circ \kappa_T\circ \pi_1(\phi_1,f(\s_0))(g) &=
\pi_1(\theta,\s_0)( g_0 e_1 g_1 \cdots e_n  g_n ) \\ &= g_0 f(e_1)
g_1 \cdots f(e_n) g_n \\ &= g
\end{align*}
and so Equation~\eqref{e:useful} holds.  Thus $\Lambda^\theta_{T,f(T)} \circ \Lambda_{f(T)} =
1$.  By conjugating Equation~\eqref{e:useful}, we obtain
\[ \Lambda_{f(T)} \circ \Lambda^\theta_{T,f(T)} = 1\]
and conclude that $\Lambda^\theta_{T,f(T)} = \Lambda_{f(T)}^{-1}$.

To show that $L^\theta_{T,f(T)} = \tilde{L}_{f(T)}^{-1}$, let
\[ u_\s = \kappa_{f(T)}^{-1}\left\{ F\theta(\pi_\s)
\left(\pi'_{f(\s)}\right)^{-1}\right\}\]
be the elements of $\pi_1(G(\Z),f(T))$ with
respect to which $L_{T,f(T)}^{\theta}$ is defined.
Here $\pi_\s$
denotes the non-backtracking path in $T$ from $\s_0$ to $\s$, and
similarly for $\pi'_{f(\s)}$ and $f(T)$.  By definition of $\theta$,
\[F\theta(\pi_\s) =
\pi'_{f(\s)} \] hence $u_\s = 1$ for all $\s \in
V(\Y)$.  Also, for each $\alpha \in \Y$, the element $h_{i(f(\alpha))}\in\fgYT$
with respect to which $\tilde{L}_{f(T)}$ is defined is a product of
oriented edges $a^{\pm}$ with $a \in T$. Hence $h_{i(f(\alpha))} =
1$.

Applying these facts, we have, for $g \in \fgYT$ and $\alpha \in
\Y$,
\begin{align*}
\tilde{L}_{f(T)} \circ L^\theta_{T,f(T)}([g],\alpha) &=
\tilde{L}_{f(T)}([\Lambda^\theta_{T,f(T)}(g)], f(\alpha)) \\ &=
\Lambda_{f(T)}\circ\Lambda^\theta_{T,f(T)}(g) h_{i(f(\alpha))}
\overline{f(\alpha)} \\ &= g ([1],\alpha) \\ &= ([g],\alpha)
\end{align*}
and
\begin{align*}
L^\theta_{T,f(T)}\circ\tilde{L}_{f(T)} ([g],f(\alpha))&=
L^\theta_{T,f(T)}(\Lambda_{f(T)}(g)\overline{f(\alpha)}) \\ &=
L^\theta_{T,f(T)}([\Lambda_{f(T)}(g)],\alpha) \\ &=
([\Lambda^\theta_{T,f(T)}\circ\Lambda_{f(T)}(g)],f(\alpha)) \\ &=
([g],f(\alpha)).
\end{align*}
Thus $L^\theta_{T,f(T)} = \tilde{L}_{f(T)}^{-1}$.
\end{pf}

The following result corresponds to Corollary~4.5,~\cite{b1:ctgg}.

\begin{cor}[Corollary to the Main Lemma] Let $G(\Y)$ and $G'(\Y')$ be developable
complexes of groups over connected scwols $\Y$ and $\Y'$, and choose
maximal trees $T$ and $T'$ in $\Y$ and $\Y'$ respectively. Suppose
$L: \DYT \to \DYpTp$ is a morphism of scwols which is equivariant
with respect to some homomorphism of groups $\Lambda: \pi_1(G(\Y),
T) \to \pi_1(G'(\Y'), T')$.  If there is a $\s_0 \in V(\Y)$ such
that
\[L([1],\s_0) = ([1],\s_0')\]
for some $\s_0' \in V(\Y')$, then there exists a morphism $\lambda:
G(\Y) \to G'(\Y')$ of complexes of groups such that
$L=L^\lambda_{T,T'}$ and $\Lambda=\Lambda^\lambda_{T,T'}$.
\end{cor}

\begin{pf} Let the quotient scwol $\Z$, the isomorphism $f:\Y \to \Z$, the data $C_\bullet$, the complex of groups $G(\Y)_{C_\bullet}$ and the isomorphism
$\theta:G(\Y) \to G(\Z)$ be as in the statement of
Lemma~\ref{l:useful} above, and similarly for $\Z'$, $f'$,
$C'_\bullet$, $G'(\Y')_{C'_\bullet}$ and $\theta'$. Let $l : \Z \to
\Z'$ be the map of quotient scwols induced by $L$ and $\Lambda$.
 By definition of $l$, $C_\bullet$ and $C_\bullet'$, and by
the assumption on $L$, we have
\[ L(\overline{\s_0}) = \overline{l(\s_0)}\]
so we may choose $N_\bullet$ with $k_{\s_0} = 1$.  Let
\[ \mu=\mu_{C_\bullet,C'_\bullet,N_\bullet}:G(\Z)_{C_\bullet} \to G'(\Z')_{C'_\bullet} \]
be the induced morphism of complexes of groups.

Let \[ \lambda = \theta'^{-1}\circ\mu\circ\theta:G(\Y) \to G'(\Y').\]
We claim that $\Lambda = \Lambda^\lambda_{T,T'}$ and $L =
L^\lambda_{T,T'}$. By
Corollary~\ref{c:commute}, it is enough to show that \[ \Lambda =
(\Lambda^{\theta'}_{T',f'(T')})^{-1} \circ
\Lambda^{\mu}_{f(T),f'(T')} \circ \Lambda^{\theta}_{T,f(T)}\]
and
\[L = (L_{T',f'(T')}^{\theta'})^{-1} \circ L_{f(T),f'(T')}^\mu \circ
L_{T,f(T)}^{\theta}. \]
The result follows from the Main Lemma applied to $\mu$, and Lemma~\ref{l:useful}
above.\end{pf}

%*******************************************************************************************************
%*******************************************************************************************************
\subsection{Coverings and developability}\label{ss:cov_dev}
%*******************************************************************************************************
%*******************************************************************************************************

This section considers the relationship between the existence of a
covering and developability.

\begin{lem}\label{l:Ydev}  Let $G(\Y)$ and $G'(\Y')$ be complexes of groups over nonempty, connected scwols $\Y$ and $\Y'$.
Assume there is a covering $\phi:G(\Y) \to G'(\Y')$.
 If $G'(\cY')$ is developable, then $G(\cY)$ is developable.
\end{lem}

\begin{pf} Let $\iota': G'(\Y') \to FG'(\Y')$ be the natural morphism defined after Definition~\ref{d:univ_gp} in
Section~\ref{sss:fundgp}.  By Proposition~\ref{p:dev_inj_FG}, since
$G'(\Y')$ is developable, $\iota'$ is injective on the local groups.
Thus, as $\phi$ is a covering, the composite morphism $\iota' \circ
\phi: G(\cY) \rightarrow FG'(\cY')$ is injective on the local
groups. Hence, by Proposition~\ref{p:dev_inj}, the complex of groups
$G(\cY)$ is developable.
\end{pf}

We do not know if the converse to Lemma~\ref{l:Ydev} holds in
general. However, in the presence of nonpositive curvature, we have
the following partial converse to Lemma~\ref{l:Ydev}.  Recall that
an $M_\kappa$-polyhedral complex is a polyhedral complex with
$n$-dimensional cells isometric to polyhedra in the simply connected
Riemannian $n$-manifold of constant sectional curvature $\kappa$.

\begin{lem}\label{l:Xnpc} Let $\phi:G(\Y) \to G'(\Y')$ be a covering of complexes of groups, over a morphism of scwols $l:\Y \to \Y'$.
Suppose that for some $\kappa \leq 0$, $\Y$ and $\Y'$ are the scwols
associated to connected $M_\kappa$-polyhedral complexes with
finitely many isometry classes of cells $K$ and $K'$ respectively,
and that $|l|:|\cY| \rightarrow |\cY'|$ is a local isometry on each
simplex.  If $G(\cY)$ has nonpositive curvature (thus is
developable), then $G'(\cY')$ also has nonpositive curvature, thus
$G'(\cY')$ is developable. \end{lem}

\begin{pf} By Lemma~\ref{l:CAT1}, to show that $G'(\cY')$ is nonpositively curved, it suffices to show that
for each vertex $\tau'$ of $K'$, the geometric link of $\tilde\tau'$
in the local development $\st(\tilde\tau')$, with the induced
spherical structure, is CAT$(1)$.  We first show, using the
following lemma, that if $\tau'$ is a vertex of $K'$, then
$\tau'=f(\tau)$ for some vertex $\tau$ of $K$.

\begin{lem}\label{l:fonto} The nondegenerate morphism of scwols $l:\cY \rightarrow \cY'$
associated to the covering $\phi:G(\cY)\rightarrow G'(\cY')$
surjects onto the set of vertices of $\cY'$.
\end{lem}

\begin{pf} Let $\sigma \in V(\cY)$ and $l(\sigma)=\sigma'\in V(\cY')$.
{}From the definitions of nondegenerate morphism of scwols and
covering of complexes of groups, it follows that every vertex of
$\cY'$ which is incident to an edge meeting $\sigma'$ lies in the
image of $l$. Since $\cY'$ is connected, we conclude that $l$
surjects onto $V(\cY')$.
\end{pf}

Let $\tau'$ be a vertex of $K'$.  By Lemma~\ref{l:fonto}, $\tau' =
l(\tau)$ for some $\tau \in V(\cY)$. Suppose $\tau$ is not a vertex
of $K$. Then there is an $a \in E(\cY)$ such that $i(a) = \tau$.  It
follows that $i(l(a))=l(i(a))=\tau'$, so $l(a) \in E(\cY')$ has
initial vertex $\tau'$. This contradicts $\tau'$ a vertex of $K'$.
Hence $\tau$ is a vertex of $K$.

Since $G(\cY)$ is nonpositively curved, the geometric link of
$\tilde\tau$ in the local development $\st(\tilde\tau)$, with the
induced spherical structure, is CAT$(1)$.  By Lemma~\ref{l:star},
there is a $\phi _\tau$-equivariant bijection $\st(\tilde\tau)
\rightarrow \st(\tilde\tau')$.  We claim this bijection is an
isometry in the induced metric, which completes the proof.

By definition of the induced metric, the action of $G_\tau$ on
$\st(\tilde\tau)$ induces a simplicial map $\st(\tilde\tau)
\rightarrow \st(\tau)$ which is a local isometry on each simplex.
Similarly, the action of $G_{\tau'}$ on $\st(\tilde\tau')$ induces
$\st(\tilde\tau) \rightarrow \st(\tau)$ which is a local isometry on
each simplex.  By assumption, the restriction of $|l|$ to
$\st(\tau)$ is a local isometry on each simplex. Hence, the
bijection $\st(\tilde\tau)\rightarrow\st(\tilde\tau')$ is a local
isometry on each simplex, and thus an isometry.
\end{pf}

%***************************************************************************************************************
%***************************************************************************************************************
%***************************************************************************************************************
\section{The Conjugacy Theorem for Complexes of Groups}\label{s:conjugacy}
%***************************************************************************************************************
%***************************************************************************************************************
%***************************************************************************************************************

In this section, we prove the analogue for complexes of groups of the Conjugacy
Theorem for graphs of groups (Theorem 5.2 of~\cite{b1:ctgg}).  Let us prove the
following lemma which characterizes coverings.

\begin{lem}[Corollary 4.6,~\cite{b1:ctgg}]\label{c:induced_covering} With the notation
in Definition~\ref{d:induced_morphism}, the induced morphism
$\lambda$ is a covering if and only if $\Lambda$ is a monomorphism
and $L$ is an isomorphism.
\end{lem}

\begin{pf} By
the Main Lemma in Section~\ref{s:covering_theory}, since the
vertical maps are isomorphisms, $\Lambda$ is a monomorphism if and
only if $\Lambda^\lambda_{T,T'}$ is a monomorphism, and $L$ an
isomorphism if and only if $L^\lambda_{T,T'}$ is an isomorphism.

Suppose $\lambda$ is a covering.  Then by
Proposition~\ref{p:induced}, $\Lambda^\lambda_{T,T'}$ is a
monomorphism and $L^\lambda_{T,T'}$ is an isomorphism, and the
conclusion follows.

Conversely, suppose $\Lambda$ is a monomorphism and $L$ is an
isomorphism.  Assume by contradiction that $\lambda$ is not a
covering.   Then either
\begin{enumerate}
\item\label{i:local2} there is a homomorphism $\lambda_\s : G_\s \to G'_{l(\s)}$ which
is not injective, or
\item\label{i:star2} there exists $a' \in
E(\cY')$ and $\s \in V(\cY)$ with $t(a')=\s'=l(\s)$, such that the
map
\[
\coprod_{\substack{a \in l^{-1}(a') \\ t(a)=\s}} G_\s
/\psi_{a}(G_{i(a)}) \to G'_{\s'}/\psi_{a'} (G'_{i(a')})\] induced by
\[
 g \mapsto \lambda_{\s}(g) \lambda(a)
 \] is not bijective.
\end{enumerate}

Condition~\eqref{i:local2} implies that the map
$\Lambda^\lambda_{T,T'}$ is not a monomorphism at $G_\s$, thus
$\Lambda^\lambda_{T,T'}$ is not a monomorphism.
Condition~\eqref{i:star2} implies that $L_{T,T'}^\lambda$ is not a
local bijection at $\St(\tilde\s)$ (see Remark~5.3,~\cite{BH}), thus
the map $L_{T,T'}^\lambda$ is not an isomorphism.  By contradiction,
we conclude that $\lambda$ is a covering.
\end{pf}

Let $\cX$ be the scwol associated to a polyhedral complex $K$. Let $G =
\Aut(K)$, and let $H$ be a subgroup of $G$ acting without inversions.  Then $H$
acts on $\cX$ in the sense of Definition~\ref{d:actiononscwol}.  Define \[G_H =
\{ g \in G \mid g\sigma \in H\sigma \mbox{ for all } \sigma \in V(\cX)\}.\] Then
$G_H$ is a subgroup of $\Aut(K)$, $H$ is a subgroup of $G_H$ and \[ H \bs \cX
= G_H \bs \cX.\]  The following theorem is the same as the Conjugacy Theorem
stated in the Introduction.

\begin{thm}\label{t:conjugacy} If $\G \leq G_H$ acts freely on $\cX$ then there is an element
$g \in G_H$ such that $g \G g^{-1} \leq H$.
\end{thm}

\begin{pf} Let $\A = H \bs \cX = G_H \bs \cX$ and $\B = \G \bs \cX$ and
let $f : \B \to \A$ be the natural projection (coming from $\G \leq
G_H$).  We form quotient complexes of groups $G(\A)=(G_\sigma,
\psi_a, g_{a,b})$ induced by the action of $H$ and
$G'(\A)=(G'_\sigma, \psi'_a, g'_{a,b})$ induced by the action of
$G_H$, using the same maximal tree $T_\A$ in the one-skeleton of
$\A$ and the same family of elements $h_a \in H \leq G_H$, for $a
\in E(\A)$. Then for each $\sigma \in V(\A)$ we have $G_\sigma \leq
G'_\sigma$, and for each edge $a \in E(\A)$ we have
$\psi_a'|_{G_{i(a)}} = \psi_a$.  There is then by
Corollary~\ref{c:induced_covering} a covering morphism
\[\lambda:G(\A) \to G'(\A)\]
induced by the identity map $L: \A \to \A$ and the inclusion
$\Lambda:H \to G_H$.
  By definition of induced
morphism, each $\lambda_\sigma$ is inclusion, and each $\lambda(a)$
is trivial.
 Hence for all $a \in E(\A)$, the inclusion induced map
\[ G_{t(a)} / \psi_a(G_{i(a)}) \to G'_{t(a)} / \psi'_a(G'_{i(a)}) \]
given by $[g] \mapsto [g]$ is a bijection.

Now form the quotient complex of (trivial) groups (since $\G$ acts
freely), $G''(\B) = (G''_\sigma, \psi''_a, g''_{a,b})$ and let
$\phi=(\phi_\sigma,\phi(b)):G''(\B) \to G'(\A)$ be a covering
morphism induced by the inclusion $\G \leq G_H$ over the natural
projection $f:\B\to \A$.  For each $a \in E(\A)$ with $t(a) =
f(\tau)$, there is a bijection
\[ \{ b \in f^{-1}(a) \mid t(b) = \tau \} \to G'_{t(a)} /
\psi_a'(G'_{i(a)})\] given by $b \mapsto [\phi(b)]$ where since
$G''(\B)$ is a complex of trivial groups, we replace the one-element
sets $G''_\tau / \psi''_b(G''_{i(b)})$ by $\{b\}$.  Thus since
$\lambda$ is a covering, for each $b \in E(\B')$, with $f(b) = a$,
we can find elements $g_b \in G_{t(a)}$ such that
\[ [\phi(b)] = [g_b] \]
in $G'_{t(a)}/\psi_a'(G'_{i(a)})$.

We now define a morphism $\phi':G''(\B) \to G(\A)$ over $f$, by each
$\phi'_\sigma$ being the inclusion of the trivial group, and
$\phi'(b) \in G_{t(f(b))}$ being $\phi'(b) = g_b$.  We then have a
bijection
\[ \{ b \in f^{-1}(a) \mid t(b) = \tau \} \to G_{t(a)} /
\psi_a(G_{i(a)})\] given by $b \mapsto [g_b]$ hence $\phi'$ is a
covering.

Choose a maximal tree $T_\B$ in $\B$ and recall that we chose a maximal tree
$T_\A$ in $\A$.  By Proposition~\ref{p:induced} the covering $\phi':G''(\B) \to
G(\A)$ induces a monomorphism of groups
\[\Lambda^{\phi'}_{T_\B,T_\A}:\pi_1(G''(\B),T_\B) \to \pi_1(G(\A),T_\A)\] and a
$\Lambda_{T_\B,T_\A}$--equivariant isomorphism of scwols
\[L^{\phi'}_{T_\B,T_\A}:D(G''(\B),T_\B) \to D(G(\A),T_\A)\] such that the
following diagram commutes: $$\xymatrix{ D(G''(\B),T_\B)
\ar[rr]^-{L^{\phi'}_{T_\B,T_\A}} \ar[d] & & D(G(\A),T_\A) \ar[d]\\ \B
\ar[rr]^-{f} & &\A}$$ where the vertical arrows are the natural projections. 
Let $\widetilde{L}_{T_\B}:D(G''(\B),T_\B) \to \cX$ and
$\widetilde{L}_{T_\A}:D(G(\A),T_\A) \to \cX$ be the canonical isomorphisms,
equivariant with respect to the isomorphisms of groups
$\Lambda_{T_\B}:\pi_1(G''(\B),T_\B) \to \G$ and
$\Lambda_{T_\A}:\pi_1(G(\A),T_\A) \to H$, respectively. Let $g \in \Aut(\cX)$ be
the following composition of isomorphisms \[ g = \widetilde{L}_{T_\A} \circ
L^{\phi'}_{T_\B,T_\A} \circ \widetilde{L}_{T_\B}^{-1}: \cX \to \cX.\] Then $g$
is equivariant with respect to the monomorphism $\theta:\G \to H$ given by \[
\theta = \Lambda_{T_\A} \circ \Lambda^{\phi'}_{T_\B,T_\A} \circ
\Lambda^{-1}_{T_\B}: \G \to H.\] Thus for all $\gamma \in \G \leq \Aut(\cX)$, we
have \[g \circ \gamma = \theta(\gamma) \circ g\] and so $g \circ \gamma \circ
g^{-1} = \theta(\gamma) \leq H$.  That is, $g \G g^{-1} \leq H$.

It remains to show that $g \in G_H$.  Let $p:\cX \to \A=H \bs
\cX=G_H \bs \cX$ and $p_\Gamma: \cX \to \B = \G \bs \cX$ be the
natural projections. Then $p = f \circ p_\G$.  We wish to show that
$p \circ g = p$.  Now $p_\G$ is the composition of
$\widetilde{L}_{T_\B}^{-1}:X \to D(G''(\B),T_\B)$ with the natural
projection $D(G''(\B),T_\B) \to \B$, and similarly $p$ is the
composition of $\widetilde{L}_{T_\A}^{-1}:\cX \to D(G(\A),T_\A)$
with the natural projection $D(G(\A),T_\A) \to \A$. Hence the
definition of $g$ and the commutativity of the diagram above mean
that $p \circ g = f \circ p_\G = p$ as required.

\end{pf}

%********************************************************************************************************
%********************************************************************************************************
%********************************************************************************************************
\section{Coverings and overgroups}\label{s:bijection}
%********************************************************************************************************
%********************************************************************************************************
%********************************************************************************************************

In this section we prove Theorem~\ref{t:bijection}, stated in the
Introduction.  We first define isomorphism of coverings.  In
Section~\ref{ss:map_a} we define a map from overgroups to coverings,
and in Section~\ref{ss:map_b} a map from coverings to overgroups.
Then in Section~\ref{ss:proof_bij} we conclude the proof of
Theorem~\ref{t:bijection} by showing that these maps are mutual
inverses.

\begin{defn}[isomorphism of coverings]\label{d:isom_covering}
Let $\lambda: G(\Y) \to G'(\Y')$ and $\lambda': G(\Y) \to G''(\Y'')$
be coverings of developable complexes of groups over connected
scwols. Fix $\s_0 \in V(\Y)$.  We say that $\lambda$ and $\lambda'$
are isomorphic coverings if for any choice of maximal trees $T$,
$T'$ and $T''$ in $\Y$, $\Y'$ and $\Y''$ respectively, there exists
an isomorphism $\lambda'': G'(\Y') \to G''(\Y'')$ of complexes of
groups such that the following diagram of morphisms of universal
covers (defined in Proposition~\ref{p:functor_morphism}) commutes
$$\xymatrix{
\DYT \ar[r]^-{L^\lambda_{T,T'}} \ar[dr]_-{L^{\lambda'}_{T,T''}} &
\DYpTp
\ar[d]^-{L^{\lambda''}_{T',T''}}\\
& \DYppTpp.}$$
\end{defn}

Note that by Corollary~\ref{c:commute}, this diagram commutes for
one triple $(T,T',T'')$ if and only if it commutes for all triples
$(T,T',T'')$.  By Proposition~\ref{p:induced}, since $\lambda$ and
$\lambda'$ are coverings, $L^\lambda_{T,T'}$ and
$L^{\lambda'}_{T,T''}$ are isomorphisms.  By
Proposition~\ref{p:functor_isom}, since $\lambda''$ is an
isomorphism, the map $L^{\lambda'\!'}_{T',T''}$ is an isomorphism.
Hence, two coverings are isomorphic if and only if they induce a
commutative diagram of isomorphisms of universal covers.

For the remainder of Section~\ref{s:bijection}, we fix the following
data:
\begin{itemize}
\item $\X$, the scwol associated to a simply connected
polyhedral complex $K$, \item $\G$, a subgroup of $\Aut(K)$
which acts on $\X$, with quotient $\Y = \G \bs \X$, \item a vertex
$\s_0 \in V(\Y)$, and
\item a set of choices $C_\bullet = (\overline\s,h_a)$ giving rise
to a complex of groups $G(\Y)_{C_\bullet}=(G_\s,\psi_a,g_{a,b})$
induced by the action of $\G$ on $\X$.
\end{itemize}

Let $\Ovg$ be the set of overgroups of $\G$ which act without
inversions, that is, the set of subgroups of $\Aut(K)$ containing
$\G$ which act without inversions.  Let  $\Covg$ be the set of
isomorphism classes of coverings of faithful, developable complexes
of groups by $G(\Y)$.

%***************************************************************************************
\subsection{The map from overgroups to coverings}\label{ss:map_a}
%***************************************************************************************

In this section we construct a map \[\underline{a}:\Ovg \to \Covg.\] We first
show in Lemma~\ref{l:o_to_c} that an overgroup induces a covering of complexes
of groups.  Then in Lemma~\ref{l:apply_main_lemma} we show that, without loss of
generality, we may apply the Main Lemma to this covering.  In
Lemma~\ref{l:a_welldefined}, we define $\underline{a}$ and show that $\underline{a}$
is well-defined on isomorphism classes of coverings.

\begin{lem}\label{l:o_to_c} Let $\Gamma'$ be an overgroup of $\Gamma$ acting without inversions.  Let $G'(\Y')_{C'_\bullet}$ be a complex of
groups over $\Y'=\G'\bs \X$ induced by the action of $\G'$ on $\X$,
for some choices $C'_\bullet$.  Let $L=Id: \X \to \X$ and let
$\Lambda: \G \hookrightarrow \Gamma'$ be inclusion, inducing $l:\Y
\to \Y'$. For some choices $N_\bullet$, let
$$\lambda=\lambda_{C_{\bullet}, C'_{\bullet}, N_{\bullet}}: G(\Y)_{C_\bullet} \to
G'(\Y')_{C'_\bullet}$$ be the morphism of complexes of groups over
$l$ induced by $L$ and $\Lambda$ (see
Definition~\ref{d:induced_morphism}).
 Then $\lambda$ is a covering.
\end{lem}

\begin{pf} By definition, $\lambda_\sigma = \Ad (k_\sigma)$, where
$k_\sigma: \overline{\sigma} \mapsto \overline{l(\sigma)}$. The
local maps $\lambda_\s$ are thus injective.

We write $[g]_a$ for the coset of $g \in G_{t(a)}$ in
$G_{t(a)}/\psi_a(G_{i(a)})$, and similarly for $[g']_{a'}$ when $g'
\in G'_{t(a')}$.  It now suffices to show that for every $a' \in
E(\Y')$ with $t(a')=\s'=l(\s) \in V(\Y)$, the map on cosets
\begin{align*}
 \coprod_{\substack{a \in l^{-1}(a')\\ t(a)=\sigma}}
\G_{\overline\s} /h_a(\G_{\overline{i(a)}})h^{-1}_a \longrightarrow
&\; \G'_{\overline{l(\s)}} /h'_{a'}
 (\G'_{\overline{i(a')}})h'^{-1}_{a'}\\
      [g]_a \longmapsto &\; [\lambda_\s (g) \lambda(a)]_{a'}
\end{align*}
is bijective. Suppose $[\lambda_\s(g) \lambda(a)]_{a'} =
[\lambda_\s(h) \lambda(b)]_{a'}$.  Then by definition of $\lambda$,
$$ h'_{a'} k_{i(b)} h_b ^{-1} k_{t(b)}^{-1} k_\sigma h^{-1}gk_\sigma ^{-1}k_{t(a)} h_a k_{i(a)} ^{-1} (h'_{a'})
^{-1} \in h'_{a'} G'_{i(a')} (h'_{a'}) ^{-1}$$ hence
$$ k_{i(b)} h_b^{-1} h^{-1}g h_a k_{i(a)}^{-1} \in \Gamma'_{\overline{i(a')}}.$$
Since $k_{i(a)}$ and $k_{i(b)}$ send $\overline{i(a)}$ and
$\overline{i(b)}$ respectively to $\overline{i(a')}$, the element
$h_b ^{-1} h^{-1}g h_a $ in $\Gamma$ sends $\overline{i(a)}$ to
$\overline{i(b)}$.  Since $l(a)=l(b)$, this implies that $a=b$.
Hence $h^{-1}g$ maps $\overline{i(a)}$ to itself, thus
$[h]_a=[g]_a$. Therefore the map on cosets is injective.

Let us show that the map on cosets is surjective. Let $[h']_{a'}$ be
an element of the target set. Let $b'=k^{-1}_\s h'
h'_{a'}(\overline{a'})$. Since $h' \in \G'_{\overline{\s'}}$, we
have $t(b')=\overline{\s}$. Let $c=p(b')$, where $p$ is the natural
projection $\X \to \Y = \G \bs \X$. Let $g \in \G_{\overline{\s}}$
be such that $g(h_c \overline{c})=b'$. We claim that $[g]_c$ maps to
$[h']_{a'}$, that is,
$$ h'^{-1}k_\s g h_c k^{-1}_{i(c)} h'^{-1}_{a'} \in h'_{a'}
\G'_{\overline{i(a')}} h'^{-1}_{a'}.$$
Since $k^{-1}_{i(c)}$ sends $\overline{i(a')}$ to $\overline{i(c)}$,
and the element $k_\s g h_c$ sends $\overline{i(c)}$ to $i(k_\s
b')=i(h'h'_{a'}(\overline{a'}))$, it follows that
$h'^{-1}_{a'}h'^{-1}k_\s g h_c k^{-1}_{i(c)}$ fixes
$\overline{i(a')}$, which proves the claim.
 \end{pf}

We now show that every covering $\lambda$ induced by an
overgroup, as in Lemma~\ref{l:o_to_c},
 is isomorphic to a covering $\lambda'$ to which the Main
Lemma may be applied.  More precisely:

\begin{lem}\label{l:apply_main_lemma} With the notation of
Lemma~\ref{l:o_to_c}, fix a vertex $\s_0 \in V(\Y)$.  Then there is a choice
$C''_\bullet$ of data for $\G'$ acting on $\X$ such that $\overline{\s_0} =
\overline{l(\s_0)}$, and a choice $N'_\bullet=\{k'_\s\}$ such that
$k'_{\s_0}=1$, so that $\lambda$ is isomorphic to the covering \[\lambda'=
\lambda'_{C_\bullet,C''_\bullet,N'_\bullet}:G(\Y)_{C_\bullet} \to
G''(\Y'')_{C''_\bullet}\] where $G''(\Y'')_{C''_\bullet}$ is the complex of
groups induced by $C''_\bullet$. \end{lem}

\begin{pf} By definition of $l$, there is a choice $C''_\bullet$ so that
$\overline{\s_0}$, determined by $C_\bullet$, equals
$\overline{l(\s_0)}$ determined by $C''_\bullet$. We now define
 a collection $N'_\bullet = \{ k'_\s\}$ such that $k'_\s
 \overline{\s} = \overline{l(\s)}$ for all $\s \in V(\Y)$.

Choose a section $s:V(\Y') \to V(\Y)$ for $l$.  That is, for each
$\s' \in V(\Y')$, choose $s(\s') \in V(\Y)$ such that $l(s(\s')) =
\s'$.  In particular, if $\s_0' = l(\s_0)$, let $s(\s_0') = \s_0$.

For each $s(\s') \in V(\Y)$, choose an element $k'_{s(\s')} \in \G'$ such that
$k'_{s(\s')} \overline{s(\s')} = \overline{\s'}$, where $\overline{s(\s')}$ is
determined by $C_\bullet$ and $\overline{\s'}$ by $C''_\bullet$.  Since
$s(\s_0') = \s_0$, and by choice of $C''_\bullet$, we have
$k'_{\s_0}\overline{\s_0} = \overline{l(\s_0)}=\overline{\s_0} $, so we may
choose $k'_{\s_0} = 1$.  For all other $\s \in V(\Y)$, let
\begin{equation}\label{e:N} k'_\s = k'_{s(l(\s))} k^{-1}_{s(l(\s))} k_\s
\end{equation} where $N_\bullet = \{ k_\s \}$.  Note that \[k'_\s \overline{\s}
= k'_{s(l(\s))} k^{-1}_{s(l(\s))} k_\s \overline{\s} = k'_{s(l(\s))}
k^{-1}_{s(l(\s))}\overline{l(\s)} = k'_{s(l(\s))}\overline{s(l(\s))} =
\overline{l(s)}.\] This defines a collection $N'_\bullet = \{ k'_\s \}$ with
$k'_{\s_0} = 1$.  Let $\lambda':G(\Y)_{C_\bullet} \to G''(\Y'')_{C'_\bullet}$ be
the covering induced by $N'_\bullet$.

We now construct an isomorphism of complexes of groups
$\mu:G'(\Y')\to G''(\Y'')$ such that the following diagram commutes
\begin{equation}\label{e:use_main}
\xymatrix{ G(\Y) \ar[r]^-{\lambda} \ar[dr]_-{\lambda'} & G'(\Y')
\ar[d]^-{\mu}\\
& G''(\Y'').}
\end{equation}
By Corollary~\ref{c:commute}, it follows that $\lambda$ is
isomorphic to $\lambda'$.

Let $f:\Y' \to \Y''$ be the identity map (both $\Y'$ and $\Y''$ are
the quotient $\G' \bs \X$).  We choose a collection $N''_\bullet =
\{ k''_{\s'} \}$ of elements of $\G'$ such that
$k''_{\s'}\overline{\s'} = \overline{f(\s')}$ as follows.  By
Equation~\eqref{e:N}, if $l(\s_1) = l(\s_2)$ then $k'_{\s_1}
k_{\s_1}^{-1} = k'_{\s_2} k_{\s_2}^{-1}$. Given $\s' \in V(\Y')$, it
is thus well-defined to put
\[k''_{\s'} = k'_\s k_\s^{-1}\]
for any $\s \in l^{-1}(\s')$.  We check
\[ k''_{\s'}\overline{\s'} = k'_\s k_\s^{-1}\overline{\s'} = k'_\s \overline{\s} = \overline{\s'} = \overline{f(\s')}\]
as required.  Define $\mu
=\mu_{C'_\bullet,C''_\bullet,N''_\bullet}:G'(\Y')_{C'_\bullet} \to
G''(\Y'')_{C''_\bullet}$.  Since $G'(\Y')$ and $G''(\Y'')$ are both
associated to the action of $\G'$ on $\X$, $\mu$ is an isomorphism.

By definition of composition of morphisms, for $g \in G_\s$ we have
\begin{align*}(\mu \circ \lambda)_\s(g) &= \mu_{l(\s)} \circ \lambda_\s (g)\\&=
\Ad(k''_{l(\s)}) \circ \Ad(k_\s) (g) \\&= \Ad(k''_{l(\s)}k_\s)(g)
\\&= \Ad (k'_\s)(g)\\& = \lambda'_\s(g)\end{align*}
and for $a \in E(\Y)$
\begin{align*}
(\mu \circ \lambda)(a) &= \mu_{l(t(a))}(\lambda(a))\mu(l(a)) \\
&=\Ad(k_{l(t(a))}'')(k_{t(a)}h_ak_{i(a)}^{-1}h_{l(a)}^{-1})k_{t(l(a))}''h_{l(a)}(k_{i(l(a))}'')^{-1}h_{f(l(a))}^{-1}\\&=
k_{l(t(a))}''k_{t(a)} h_a k_{i(a)}^{-1}
(k_{i(l(a))}'')^{-1}h_{f(l(a))}^{-1}\\& = k'_{t(a)} h_a
(k_{i(a)}')^{-1} h_{f(l(a))} \\&= \lambda'(a)
\end{align*}
hence the diagram at~\eqref{e:use_main} commutes.
\end{pf}

\begin{lem}\label{l:a_welldefined} Let
\[\underline{a}:\Ovg \to \Covg\] be a map taking an overgroup $\G'$ of
$\G$ to a covering, as described in Lemma~\ref{l:o_to_c}.
Let $C'_\bullet$, $N_\bullet$ and $C''_
\bullet$, $N'_\bullet$ be any two choices for the construction of
$\underline{a}(\G')$
\[\lambda_{C_\bullet,C'_\bullet,N_\bullet}:G(\Y)_{C_\bullet} \to
G'(\Y')_{C'_\bullet} \quad \mbox{and}\quad
\lambda'_{C_\bullet,C''_\bullet,N'_\bullet}:G(\Y)_{C_\bullet} \to
G''(\Y'')_{C''_\bullet}\] Then
$\lambda$ and $\lambda'$ are isomorphic coverings, so $\underline{a}$ is
well-defined.
\end{lem}

\begin{pf}  Fix a vertex $\s_0 \in V(\Y)$ and let $\s_0'=l(\s_0)$. By
Lemma~\ref{l:apply_main_lemma}, we may without loss of generality assume that
the Main Lemma may be applied to $\lambda$ and $\lambda'$.  As in the proof of
Lemma~\ref{l:apply_main_lemma}, choose a collection $N''_\bullet=\{k''_{\s'}\}$
with $k''_{\s_0'} = k'_{\s_0}k_{\s_0}^{-1} =1$.  Then we may apply the Main
Lemma to the isomorphism of complexes of groups \[
\lambda''=\lambda_{C_\bullet',C_\bullet'',N_\bullet''}:G'(\Y')_{C_\bullet'} \to
G''(\Y'')_{C''_\bullet}.\]

Choose maximal trees $T$, $T'$ and $T''$ in $\Y$, $\Y'$ and $\Y''$
respectively.  We need to check that the triangle
\begin{equation}\label{E:triangle}
\xymatrix{ \DYT \ar[r]^-{L^\lambda_{T,T'}}
\ar[dr]_-{L^{\lambda'}_{T,T''}} & \DYpTp
\ar[d]^-{L^{\lambda''}_{T',T''}}\\
& \DYppTpp}
\end{equation}
commutes.  Using the Main Lemma three times, we obtain the diagram
$$\xymatrix{
\DYT \ar[dd]_{\tilde{L}_{T}} \ar[drrr]_{L^{\lambda'}_{T,T''}}
\ar[rr]^{L^\lambda _{T,T''}} & & \DYpTp \ar[dd]^{\tilde{L}_{T'}}
\ar[dr]^{L^{\lambda''}_{T',T''}}&\\
&&& \DYppTpp \ar[dd]_{\tilde{L}_{T''}}\\
\cX \ar[rr]^{Id} \ar[drrr]^{Id} &&\cX \ar[dr]^{Id}\\
&&& \cX, \\}
$$ and see that the commutativity of (\ref{E:triangle}) is
equivalent to the commutativity of the tautological triangle
$$\xymatrix{
\cX \ar[r]^-{Id} \ar[dr]_-{Id} & \cX \ar[d]^-{Id}\\
& \cX}$$ which is obvious. \end{pf}

%**********************************************************************
\subsection{The map from coverings to overgroups}\label{ss:map_b}
%**********************************************************************

We now show that there is a map \[\underline{b}:\Covg \to \Ovg\] Let
$\lambda: G(\Y) \to G'(\Y')$ be a covering of complexes of groups,
where $G'(\Y')$ is faithful and developable.  For any maximal
subtrees $T$ and $T'$ of $\Y$ and $\Y'$ respectively, let
$\Lambda_{T,T'} : \pi_1(G(\Y),T) \to \pi_1(G'(\Y'), T')$ be the
associated group monomorphism, and $L^\lambda_{T,T'}: \DYT \to
\DYpTp$ be the associated $\Lambda_{T,T'}$-equivariant isomorphism
of scwols. Composition with the isomorphism $\tilde{L}_{T}^{-1}$
(see Proposition~\ref{p:isoms}) yields an isomorphism of scwols
$$L_{\lambda, T'}= L^\lambda_{T,T'} \circ \tilde{L}_{T}^{-1} : \mathcal{X} \to \DYpTp$$
which is equivariant with respect to $\Lambda_{T,T'} \circ
\Lambda_{T}^{-1}: \G \to \pi_1(G'(\Y'),T')$.  We set
$\underline{b}(\lambda)$ to be the group
\[\underline{b}(\lambda)=L_{\lambda,T'}^{-1}(\pi_1(G'(\Y'),T'))L_{\lambda,T'}\]
which acts on $\X$. Since $G'(\Y')$ is faithful, $\fgYpTp$ acts
faithfully on $\DYpTp$.  Hence we may identify
$\underline{b}(\lambda)$ with a subgroup of $\Aut(K)$ which acts on
$\X$. As $\Lambda_{T,T'}$ is injective, $\underline{b}(\lambda)$ is
an overgroup of $\G$.

Lemma~\ref{l:b_welldefined} below shows that $\underline{b}$ is
well-defined, that is, only depends on the isomorphism class of the
covering $\lambda$.

\begin{lem}\label{l:b_welldefined} Let $\lambda: G(\Y) \to
G'(\Y')$ and $\lambda':G(\Y) \to G''(\Y'')$ be isomorphic coverings
of complexes of finite groups, with $G'(\Y')$  and $G''(\Y'')$
faithful and developable.  Then
$\underline{b}(\lambda)=\underline{b}(\lambda')$.\end{lem}

\begin{pf}
By definition, there exists an isomorphism $\lambda'': G'(\Y') \to
G''(\Y'')$ such that, for any choice of maximal trees, we have a
commuting triangle
$$\xymatrix{
\DYT \ar[r]^{L^\lambda_{T,T'}} \ar[dr]_-{L^{\lambda'}_{T,T''}} &
\DYpTp
\ar[d]^-{L^{\lambda''}_{T',T''}}\\
& \DYppTpp}$$ and thus, composing with $\tilde{L}_{T}^{-1}$, a
commuting triangle
$$\xymatrix{
\mathcal{X} \ar[r]^{L_{\lambda,T'}\;\;} \ar[dr]_-{L_{\lambda',T''}}
&
\DYpTp \ar[d]^-{L^{\lambda''}_{T',T''}}\\
& \DYppTpp.}$$  Since $\lambda''$ is an isomorphism, by
Proposition~\ref{p:functor_isom} the group homomorphism $
\Lambda_{T',T''}:\fgYpTp \to \fgYppTpp $ is an isomorphism.  Thus,
as $L^{\lambda''}_{T',T''}$ is $\Lambda_{T',T''}$-equivariant,
\begin{align*}
\underline{b}(\lambda') &= L_{\lambda',T''}^{-1}(\pi_1(G''(\Y''),
T''))L_{\lambda',T''}^{-1}\\&=L_{\lambda,T'}^{-1}
(L^{\lambda''}_{T',T''})^{-1}(\pi_1(G''(\Y''),
T''))L^{\lambda''}_{T',T''}L_{\lambda,T'}
\\&=L_{\lambda,T'}^{-1}(\pi_1(G'(\Y'),
T'))L_{\lambda,T'} \\&=\underline{b}(\lambda).\end{align*} Therefore
$\underline{b}$ is well-defined.\end{pf}

%********************************************************************
%********************************************************************
\subsection{Proof of Theorem~\ref{t:bijection}}\label{ss:proof_bij}
%********************************************************************
%********************************************************************

We now complete the proof of Theorem~\ref{t:bijection}.  Let
$\underline{a}:\Ovg \to \Covg$ be as defined in
Section~\ref{ss:map_a} and $\underline{b}:\Covg \to \Ovg$ be as
defined in Section~\ref{ss:map_b}.

\begin{prop}\label{p:bijection} The maps $\underline{a}$ and $\underline{b}$
are mutually inverse bijections. \end{prop}

\begin{pf} We first prove that $\underline{b}\circ \underline{a}=1$. For
this, let $\Gamma'$ be an overgroup of $\Gamma$ acting without
inversions, and let $\underline{a}(\Gamma')=\lambda: G(\Y)  \to
G'(\Y')$ be an associated covering over a morphism of scwols $l:\Y
\to \Y'$.  By
Lemma~\ref{l:apply_main_lemma}, we may assume that
we can apply the Main
Lemma to $\lambda$. For any maximal subtrees $T$ and $T'$ of
$\Y$ and $\Y'$ respectively, we have then a commuting diagram of
(equivariant) isomorphisms of scwols
$$\xymatrix{
\DYT \ar[d]_-{\tilde{L}_T} \ar[r]^-{L^\lambda_{T,T'}} & \DYpTp
\ar[d]^-{\tilde{L}_{T'}}\\
\mathcal{X} \ar[r]_-{L=Id} & \mathcal{X}.}$$ Thus
\begin{align*}
\underline{b}(\lambda) &=
L_{\lambda,T'}^{-1}(\pi_1(G'(\Y'),T'))L_{\lambda,T'}\\&=
(L^\lambda_{T,T'}\circ\tilde{L}_T^{-1})^{-1}(\pi_1(G'(\Y'),T'))L^\lambda_{T,T'}\circ\tilde{L}_T^{-1}\\
&=\tilde{L}_{T'}
(\pi_1(G'(\Y'),T'))\tilde{L}_{T'}^{-1}\\&=\Gamma'\end{align*} since
$\tilde{L}_{T'}$ is equivariant with respect to the isomorphism
$\Lambda_{T'}:\fgYpTp \to \G'$.  We conclude that $\underline{b}\;
\underline{a}(\Gamma')=\Gamma'$.

We now prove that $\underline{a} \circ \underline{b}=1$. Let
$\lambda: G(\Y) \to G'(\Y')$ be a covering of a faithful developable
complex of groups $G'(\Y')$ over a morphism of scwols $l:\Y \to
\Y'$.  Choose a vertex $\s_0 \in V(\Y)$ and
maximal trees $T$ and $T'$ in $\Y$ and $\Y'$ respectively.
Without loss of generality, we identify $G'(\Y')$ with the complex of
groups induced by the action of $\fgYpTp$ on $\DYpTp$, using the isomorphism
$\theta'$ defined in Lemma~\ref{l:useful} above.  By abuse of notation, we write
$\lambda$ for $\theta'\circ\lambda$.  Let $\G'=\underline{b}(\lambda)$.

 Let $\mu =
\underline{a}(\Gamma')$ be a covering $\mu: G(\Y) \to
G''(\Y'')_{C''_\bullet}$ over a morphism of scwols $l':\Y \to \Y''$,
where $G''(\Y'')$ is a complex of groups induced by the action of
$\G'$ on $\X$.  By Lemma~\ref{l:apply_main_lemma}, we may assume
that $\overline{\s_0}=\overline{l'(\s_0)}$ so that we can apply the
Main Lemma to $\mu$. We now show that $\lambda$ and
$\mu=\underline{a}\,\underline{b}(\lambda)$ are isomorphic
coverings.

The map $\underline{b}$ induces a group isomorphism
\[\Lambda_{\underline{b}}:\fgYpTp \to \underline{b}(\lambda)\]
with, for each $g' \in \fgYpTp$ and each $\alpha \in \X$,
\[\Lambda_{\underline{b}}(g')\cdot \alpha = L_{\lambda,T'}^{-1}( g' \cdot
L_{\lambda,T'}(\alpha)).\]
By construction, $L_{\lambda,T'}^{-1}:\DYpTp \to \X$ is
$\Lambda_{\underline{b}}$-equivariant.  Let $f:\Y' \to \Y''$ be the induced map
of the quotient scwols
\[\Y' = \fgYpTp \bs \DYpTp \quad\mbox{and}\quad \Y'' = \G' \bs \X.\]
 Since
$\Lambda_{\underline{b}}$ and $L_{\lambda,T'}^{-1}$ are both isomorphisms, $f$
is an isomorphism of scwols. We claim that the following diagram of morphisms of
scwols commutes:
$$\xymatrix{
\Y \ar[r]^{l} \ar[dr]_-{l'} &
\Y'
\ar[d]^-{f}\\
& \Y''.}$$ Let $\alpha \in \Y$.  Then $\alpha = \G \overline{\alpha}$ with
$\overline{\alpha} \in \X$.  We identify $l(\alpha) \in \Y'$ with the orbit
$\fgYpTp ([1],l(\alpha))=\fgYpTp([u_{i(\alpha)}],l(\alpha))$.  Then
\[f(l(\alpha)) = \G' L_{\lambda,T'}^{-1}([u_{i(\alpha)}], l(\alpha)) = \G'
h_{i(\alpha)} \overline\alpha = \G' \overline\alpha = l'(\alpha)\] proving the
claim.

We next choose elements $k_{\s'} \in \G'$ such that, for each $\s'
\in V(\Y')$,
\[ k_{\s'} \, L^{-1}_{\lambda,T'}([1],\s') = \overline{f(\s')}. \]
We claim that $L^{-1}_{\lambda,T'}([1],l(\s_0)) =
\overline{f(l(\s_0))}$.  Now
\[L_{\lambda,T'}(\overline{f(l(\s_0))}) = L_{T,T'}^\lambda \circ \tilde{L}_T^{-1}(\overline{f(l(\s_0))})=L_{T,T'}^\lambda([1],\s_0) = ([1],l(\s_0))\]
since $h_{i(f(l(\s_0)))}=1$ and $u_{\s_0} = 1$, which proves the
claim.  Hence we may, and do, choose $k_{\s_0'} = 1$.

The elements $k_{\s'}$ then induce a morphism
 $\phi:G'(\Y') \to G''(\Y'')$ over $f$, given by $\phi_{\s'}(g') =
 k_{\s'}\Lambda_{\underline{b}}(g')k_{\s'}^{-1}$ for $g' \in G'_{\s'}$, and
 $\phi(a') = k_{t(a')}
 \Lambda_{\underline{b}}(a'^+)k_{i(a')}^{-1}h_{f(a')}^{-1}$ for $a' \in E(\Y')$.
 Since
 $\Lambda_{\underline{b}}$ and $f$ are isomorphisms, $\phi$ is an isomorphism of
 complexes of groups.
 Moreover, the following diagram commutes up to a homotopy from
 $\Lambda_{\underline{b}} \iota'_{T'}$ to $\phi_1'' \phi$, given by the elements
 $\{k_{\s'}\}$:
$$\xymatrix{
G'(\Y') \ar[r]^-{\iota'_{T'}} \ar[d]_-{\phi} & \fgYpTp
\ar[d]_-{\Lambda_{\underline{b}}} \\
G''(\Y'') \ar[r]^-{\phi_1''} & \G'.}$$ Hence, by
Proposition~\ref{p:functor_dev}, there is a
$\Lambda_{\underline{b}}$-equivariant isomorphism of scwols
\[L_{\underline{b}}: \DYpTp \to D(\Y'',\phi_1'')\]
given explicitly by
\[([g'],\alpha') \mapsto
([\Lambda_{\underline{b}}(g')k_{i(\alpha')}^{-1}],f(\alpha')).\]

We now choose a maximal subtree $T''$ of $\Y''$ and
compose $L_{\underline{b}}$ with the isomorphism $L_{T''}^{-1}:
D(\Y'',\phi_1'') \to D(\Y'',T'')$ to obtain an isomorphism of scwols
\[L: \DYpTp \to \DYppTpp\]
which is equivariant with respect to the composition of group isomorphisms
\[\Lambda_{T''}^{-1} \circ \Lambda_{\underline{b}} : \fgYpTp \to \G' \to
\fgYppTpp.\]  Since $k_{\s_0'} = 1$ and $h_{f(\s_0')}=1$,
\[ L_{\underline{b}}([1],\s_0') = ([k_{\s_0'}], f(\s_0') ) =
([h_{f(\s_0')}],f(\s_0')) = L_{T''}([1],f(\s_0'))\] hence
$L([1],\s_0') = ([1],f(\s_0'))$.  We may thus apply the Corollary to
the Main Lemma to $L$.  We now have $L=L^{\lambda'}_{T',T''}$ for some
morphism $\lambda': G'(\Y') \to G''(\Y'')$.  By
Proposition~\ref{p:functor_isom}, since $L$ is an isomorphism of
scwols which is equivariant with respect to an isomorphism of
groups, $\lambda'$ is an isomorphism of complexes of groups.

To complete the proof, it now suffices to show that the following diagram
commutes:
$$\xymatrix{
\DYT \ar[r]^{L^\lambda_{T,T'}} \ar[dr]_-{L^{\mu}_{T,T''}} &
\DYpTp
\ar[d]^-{L=L^{\lambda'}_{T',T''}}\\
& \DYppTpp.}$$
By definition of $L$, it suffices to show that \[L_{\underline{b}} \circ
L_{T,T'}^\lambda = L_{T''} \circ L_{T,T''}^\mu.\]  Let $g \in \fgYT$ and
$\alpha \in \Y$.  We write $u^\lambda_{i(\alpha)}$ for the element of
$\fgYpTp$ with respect to which $L^\lambda_{T,T'}$ is defined, and similarly
for $u^\mu_{i(\alpha)} \in \fgYppTpp$.  Then
\[L_{\underline{b}}\circ L^\lambda_{T,T'}([g],\alpha) =
\left(\left[\Lambda_{\underline{b}}\left\{\Lambda_{T,T'}(g)
u^\lambda_{i(\alpha)}\right\}k^{-1}_{i(l(\alpha))}\right],f(l(\alpha))\right)\]
and
\[L_{T''} \circ L_{T,T''}^\mu([g],\alpha) =
\left(\left[\Lambda_{T''}\left\{\Lambda_{T,T''}(g)
u^\mu_{i(\alpha)}\right\}h_{i(l'(\alpha))}\right],l'(\alpha)\right).\]
Since $f \circ l = l'$, it suffices to show that
\begin{equation}\label{e:commute}\Lambda_{\underline{b}}\left\{\Lambda_{T,T'}(g)
u^\lambda_{i(\alpha)}\right\}k^{-1}_{i(l(\alpha))} \overline{f(l(\alpha))}
= \Lambda_{T''}\left\{\Lambda_{T,T''}(g)
u^\mu_{i(\alpha)}\right\}h_{i(l'(\alpha))}\overline{l'(\alpha)}.\end{equation}
By definition of the elements $k_{\s'}$, the left-hand side of~\eqref{e:commute} equals
\begin{align*}
\Lambda_{\underline{b}}&\left\{\Lambda_{T,T'}(g)
u^\lambda_{i(\alpha)}\right\} L_{\lambda,T'}^{-1}\left([1],l(\alpha)\right) \\& =
L_{\lambda,T'}^{-1} \left(
\Lambda_{T,T'}(g)u^\lambda_{i(\alpha)}\cdot\left([1],l(\alpha)\right) \right)
\quad\mbox{since $L_{\lambda,T'}^{-1}$
is $\Lambda_{\underline{b}}$-equivariant} \\ &= L_{\lambda,T'}^{-1}
\left( \Lambda_{T,T'}(g) \cdot( [u^\lambda_{i(\alpha)}],l(\alpha))\right)
\\ &= L_{\lambda,T'}^{-1}  \left( \Lambda_{T,T'}(g) \cdot
L^\lambda_{T,T'}([1],\alpha) \right) \\ &=
L_{\lambda,T'}^{-1}  \circ L^\lambda_{T,T'}  ([g],\alpha)
\quad\mbox{since  $L^\lambda_{T,T'}$ is $\Lambda_{T,T'}$-equivariant}\\
& = \tilde{L}_T([g],\alpha)\quad\mbox{by definition of $ L_{\lambda,T'}$.}
\end{align*}
On the right-hand side of~\eqref{e:commute}, we have, by definition
of $\tilde{L}_{T''}$,
\begin{align*}
\Lambda_{T''}&\left\{\Lambda_{T,T''}(g)
u^\mu_{i(\alpha)}\right\}\tilde{L}_{T''}([1],l'(\alpha)) \\
&= \tilde{L}_{T''}\left( \Lambda_{T,T''}(g)
u^\mu_{i(\alpha)} \cdot([1],l'(\alpha))\right)\quad\mbox{since $\tilde{L}_{T''}$ is
$\Lambda_{T''}$-equivariant}\\
&=\tilde{L}_{T''}\left( \Lambda_{T,T''}(g)\cdot
 ([u^\mu_{i(\alpha)}],l'(\alpha))\right)\\
&= \tilde{L}_{T''}\left( \Lambda_{T,T''}(g)\cdot
 L^\mu_{T,T''}([1],\alpha)\right)\\
&= \tilde{L}_{T''}\circ L^\mu_{T,T''}([g],\alpha) \quad\mbox{since
${L}^\mu_{T,T''}$ is
$\Lambda_{T,T''}$-equivariant.}
\end{align*}
But by the Main Lemma applied to $\mu$, we have a commuting
square
$$\xymatrix{
\DYT \ar[d]_-{\tilde{L}_T} \ar[r]^-{L^\mu_{T,T''}} & \DYppTpp
\ar[d]^-{\tilde{L}_{T''}}\\
\mathcal{X} \ar[r]_-{Id} & \mathcal{X}}$$
hence equation~\eqref{e:commute} holds.
 \end{pf}

We conclude by establishing a bijection between $n$--sheeted coverings and
overlattices of index $n$.

\begin{cor}\label{c:main} Let $K$ be a simply connected, locally finite
polyhedral complex, and let $\G$ be a cocompact lattice in $\Aut(K)$
(acting without inversions) which induces a complex of groups
$G(\Y)$.  Then there is a bijection between the set of overlattices
of $\Gamma$ of index $n$ (acting without inversions)  and the set of
isomorphism classes of $n$--sheeted coverings of faithful developable
complexes of groups by $G(\Y)$. \end{cor}

\begin{pf} By the definition of $n$--sheeted covering, the bijection of
Theorem~\ref{t:bijection} sends an isomorphism class of
finite-sheeted coverings to an overgroup containing $\G$ with finite
index.

Since $\G$ is cocompact, the quotient scwol $\Y$ is finite and the
local groups $G_\s$ of $G(\Y)$ are finite groups.  Let
$\lambda:G(\Y) \to G'(\Y')$ be a finite-sheeted covering, where
$G'(\Y')$ is a faithful, developable complex of groups.  Then $\Y'$
is finite by Lemma~\ref{l:fonto}, and the local groups $G'_{\s'}$
are finite since $\lambda$ is finite-sheeted. It follows that the
overgroup $\underline{b}(\lambda)$ is a cocompact lattice acting
without inversions on $K$.

It remains to show that the bijection $\underline{a}$ sends an
overlattice $\Gamma'$ of index $n$ to an $n'$-sheeted covering, with
$n=n'$.  Let $\lambda=\underline{a}(\G'):G(\Y) \to G'(\Y')$ be a
covering associated to $\G'$, over the morphism of quotient scwols
$l: \G \ba \cX \to \G' \ba \cX$.   Then
\begin{align*}
n&=[\Gamma': \Gamma]=\frac{\Vol(\Gamma \quot V(\cX))}{\Vol(\Gamma'
\quot V(\cX))}=
\frac{\sum_{ \s \in V(\Y)}\frac{1}{|G_\s|}}{\sum_{\s' \in V(\Y')}
\frac{1}{|G'_{\s'}|} } \\
&=\frac{ \sum_{ \s' \in V(\Y')} \; \sum_{ \s \in
l^{-1}(\s')} \frac{1}{|G_\s|}}{ \sum_{ \s' \in
V(\Y')}\frac{1}{|G'_{\s'}|}}= \frac{\sum_{\s' \in
V(\Y')} \frac{n'}{|G'_{\s'}|}}{\sum_{\s' \in V(\Y')}
\frac{1}{|G'_{\s'}|}}=n'
\end{align*} as required.\end{pf}

We remark that we can define isomorphism between two coverings $\lambda: G'(\Y') \to G(\Y)$
and $\lambda: G''(\Y'') \to G(\Y)$ analogous to Definition~\ref{d:isom_covering}
so that there is a bijection between the set of subgroups of $\Gamma$ of index
$n$ and the set of isomorphism classes of $n$--sheeted coverings of $G(\Y)$ by
faithful developable complexes of groups. Since the proof is similar to that of
Corollary~\ref{c:main}, we omit it. Note that the developability comes free by
Lemma~\ref{l:Ydev}.

\begin{ack}  We thank Andr\'e Haefliger and Martin Bridson for helpful
correspondence, and Lisa Carbone for encouragement.  We are grateful
to Fr\'ed\'eric Paulin and Benson Farb for their constant help, advice and
support. We also thank Yale University and the University of Chicago for mutual visits which enabled this work. \end{ack}

%*****************************************************************************

\end{document}